\documentclass[11pt]{article}
\usepackage{amsmath,amssymb,amsthm,amscd,tabularx,array}
   
\newfont{\extra}{msbm10 scaled\magstep1}

\newcommand{\extr}[1]{\mbox{\extra #1}}

\newcommand{\sect}[1]{\setcounter{equation}{0}\section{#1}}
\newcommand{\subsect}[1]{\subsection{#1}}
\newcommand{\subsubsect}[1]{\subsubsection{#1}}

\newtheorem{theorem}{Theorem}[section] 
\newtheorem{proposition}{Proposition}[section] 
 
\newtheorem{corollary}{Corollary}[section]

\def\be{\begin{equation}}
\def\ee{\end{equation}}
\def\bea{\begin{eqnarray}}
\def\eea{\end{eqnarray}}

\newcommand{\comb}[2]{\left( \begin{array}{c} #1 \\ #2 \end{array}\right)}

\newcommand{\lad}[1]{ {\text{\rm ad}^{\text{\rm l}}_{#1}}}
\newcommand{\rad}[1]{ {\text{\rm ad}^{\text{\rm r}}_{#1}}}

\def\C{\extr C}
\def\K{\extr K}
\def\N{\extr N}
\def\R{\extr R}
\def\v{\varepsilon}
\def\k{\kappa}

\newcommand{\bicross}{\triangleright\!\!\!\blacktriangleleft}
\newcommand{\rimo}{\triangleright\!\!\!<}
\newcommand{\leco}{>\!\!\blacktriangleleft}
\newcommand{\act}{\triangleright}
\newcommand{\RL}{\triangleright\!\!\!\blacktriangleleft}
\newcommand{\LR}{\blacktriangleright\!\!\!\triangleleft}
\newcommand{\LLL}{\triangleleft}

\newcommand{\RIMO}{\triangleright\!\!\!<}

\newcommand{\LEMO}{>\!\!\!\triangleleft}
\newcommand{\RICO}{\blacktriangleright\!\!\!<}

\newcommand{\lact}{\triangleright}
\newcommand{\ract}{\triangleleft}
\newcommand{\lcact}{\blacktriangleleft}
\newcommand{\rcact}{\blacktriangleright}

\begin{document}

\begin{center}{ \LARGE \bf
 Representations 
of \\[0.4cm]
Quantum  Bicrossproduct Algebras}
\end{center}
\vskip0.25cm

\begin{center}
Oscar Arratia $^1$ and Mariano A. del Olmo $^2$
\vskip0.25cm

{ \it $^{1}$ Departamento de  Matem\'atica Aplicada a la  Ingenier\'{\i}a,  \\
Universidad de  Valladolid. E-47011, Valladolid, Spain. 
\vskip0.15cm
$^{2}$ Departamento de  F\'{\i}sica Te\'orica,\\
 Universidad de  Valladolid, 
 E-47011, Valladolid,  Spain}

\vskip0.15cm

E. mail: oscarr@wmatem.eis.uva.es, olmo@fta.uva.es
 
\end{center}
 
\vskip1cm
\centerline{\today}
\vskip.75cm

\begin{abstract}
We present a method to construct 
induced representations of quantum algebras having the structure of
bicrossproduct.  We apply this procedure to some quantum kinematical algebras in $(1+1)$
dimensions with  this kind of structure: null-plane quantum  Poincar\'e algebra,
non-standard quantum  Galilei algebra and quantum kappa Galilei algebra.
\end{abstract}
\vskip .75cm

\newpage

\sect{Introduction}
\label{introduction}

In a recent paper \cite{olmo00} we developed a  method to construct
induced representations of quantum algebras mainly based on the concepts of module and
duality. Since by dualization objects like modules and comodules can be seen as
equivalent, then we have not only  regular and induced representations but also coregular
and coinduced representations. 
The main result of that work was the possibility of constructing coregular and coinduced
representations of a Hopf algebra $U_q(\mathfrak g)$  when dual bases of it and its dual
$Fun_q(G)$ (or $F_q(G)$) are known, being $\mathfrak g$ the Lie algebra of a Lie group $G$.

Now we present a particular study of the induction  method for those quantum
algebras having a structure of bicrossproduct, which
is a  generalization of the idea of semidirect product of Lie groups to the quantum
case.  
This kind of structure of semidirect product is well known in physics where many
interesting groups, like Euclidean,  Galilei and Poincar\'e have it. 
The corresponding quantum Lie algebras inherit  the `semidirect'
structure in the algebra sector and the algebra of functions also has a semidirect
product structure in the coalgebra sector.  These ideas were generalized by  Molnar
\cite{Molnar} with the notions of  smash product or, more recently by Majid with that of
bicrossproduct \cite{Maj88}--\cite{Maj90d}. 

The quantum counterparts of the above mentioned groups and algebras are related with the 
symmetries of the physical space-time in a noncommutative framework. The study of 
these quantum symmetries and their representations generalizes the well known and fruitful  programm 
started by Wigner in 1939 \cite{wigner39} inside the perspective of the  noncommutative geometry
\cite{connes}, which  in the last years is finding many applications  in 
physics (see, for instance, \cite{connes01} and references therein).

In this paper we continue the analysis of the theory of induced representations but now
referred  to Hopf algebras with bicrossproduct structure, whose first factor is
cocommutative and the second one commutative. Our fundamental objective is the
description of the representations induced by characters of the commutative sector.  We
want to avoid the problems derived by the use of pairs of dual bases and to open new ways
which allow to do some computations, whose difficulty increases with the number of
generators. However, although the first results are also obtained using dual bases, they 
show the existence of some lying structure connected with some one-parameter flows
defined by the cocommutative sector over an object related with the commutative factor.
The nature of this structure will be clear after an adequate reinterpretation of the
factors of the bicrossproduct. More explicitly, the cocommutative factor will be seen as 
the enveloping algebra of a certain Lie algebra but the commutative factor will be
identified with the algebra of functions of another Lie group. In this way,  the action,
defining a part of the bicrossproduct structure in the original Hopf algebra, is the
result of translating to the algebra of functions the action of a Lie group over another Lie
group.   

A crucial point in our approach is to describe the Hopf algebra substituting the monomial
bases by elements which are product  of an element of the group associated to the first
factor times a function  belonging to the second one. That allows to prove theorems
\ref{crbicross} and \ref{repequiv} that are the cornerstone of this paper.  The first theorem describes
the four regular modules associated to a bicrossproduct Hopf algebra in terms
of the regular actions of its components and the action, mentioned above, associated with the
bicrosproduct structure. The second theorem allows the construction of the representations induced by
characters of the abelian sector and classifies the equivalence classes of the induced representations 
in terms of the orbits associated to the action  of a certain group. Moreover, a
$*$--structure is introduced  in such a way that the induced representations are unitary.

The induction procedure, such as it has been formulated by  us, has an algebraic character
since  we  use  objects like modules, comodules, etc.,  which are
the appropriate tools to work with the algebraic structures  exhibited by the quantum
groups and algebras \cite{olmo00,olmo98,oscarTesis,ao00}.

Dobrev  has  developed \cite{dobrev1,dobrev2} a method
for construct representations of quantum groups  similar, in some sense, to ours, i.e., 
both methods  emphasize the dual case, closer to the classical one, and the
representations are constructed in the algebra sector.  
We can mention some works that  have also extended 
the induction technique to
quantum groups but constructing corepresentations, i.e. representations of the coalgebra
sector  \cite{ibort}--\cite{bgst98}.   

The paper is organized as follows. Section
\ref{mathematical preliminaries} is devoted to review the main ideas and concepts that we
will use along the paper, like module, comodule, module algebra, bicrossproduct, etc. 
The last part of this section presents original results showing how  pairs of dual bases
and $*$--structures over bicrossproduct Hopf algebras can be obtained starting from those
of their factors. The first results about induced representations of quantum algebras with
bicrossproduct structure are presented in  section \ref{induced representations}. We
obtain the representations  making use of pairs of dual bases. In section
\ref{modulesandrepresentations} we begin to study the induction problem  taking into
account the deep relation between modules and representations obtaining, in some sense,
more deep results from a geometric point of view using the concept of regular co-space.
In section \ref{ejemplos} we obtain the induced representations of some kinematical
quantum algebras making use of the method developed in the previous sections. We end with
some comments and conclusions.

\sect{Mathematical preliminaries}
\label{mathematical preliminaries}

Let $H=\left({V};{m}\ {\eta};{\Delta}\
{\epsilon};{S}\right)$ be a Hopf algebra with underlying vector space $V$ over the field
$\K$ ($\C$ or $\R$), multiplication 
$m: H \otimes H \to H$, coproduct $\Delta : H \to H\otimes H$, unit $ \eta:
\K \to H $, counit $ \epsilon: H \to \K$ and  antipode $S: H \to
H$. 

A  Hopf algebra  can  be considered as a bialgebra with an antilinear map $S$, and a
bialgebra can be seen as composed by two `substructures' or `sectors' (the algebra sector 
$(V, m, \eta)$ and the coalgebra sector $(V,\Delta, \epsilon)$) with some compatibility
conditions \cite{ChP}).

On the other hand, the algebras considered in this work are finitely generated although
they are infinite dimensional.  For this reason the following  multi-index notation is
very useful \cite{olmo00}. 
Let us suppose that $A$ is an algebra
generated by the  elements $(a_1, a_2, \ldots, a_r)$
and the ordered monomials 
\begin{equation} 
a^n := a_1^{n_1} a_2^{n_2} \cdots a_r^{n_r} \in A, \qquad 
n=(n_1, n_2, \ldots, n_r) \in \mathbb{N}^r,
\end{equation}
form a basis of the linear space  underlying to $A$.
An arbitrary product of generators of $A$ is written in a normal ordering if it is 
expressed in terms of the basis  $(a^n)_{n \in \mathbb{N}^r}$. In some
cases we will use the notation $(a_n := a_1^{n_1} a_2^{n_2} \cdots a_r^{n_r})_{n \in
\mathbb{N}^r}$. For
$0=(0,\ldots, 0) \in
\N^n$ we have $a^0\equiv 1_A$. Multi-factorials
and multi-deltas are defined by
 \begin{equation}
  l! = \prod_{i=1}^n l_i!, \qquad \delta_l^m= \prod_{i=1}^n \delta_{l_i}^{m_i}.
\end{equation}

\subsect{Duality}
\label{duality}

It is well known that the dual object of $V$ is defined
as the vector space of its linear forms, i.e., 
$V^*= {\cal L}(V, \K)$. Hence, if $(V,m,\eta)$ is a finite algebra it is natural to
define the dual objet as $(V^*, m^*, \eta^*)$ obtaining 
a coalgebra and viceversa. However, in the infinite dimensional case the spaces $(V
\otimes V)^*$ and $V^* \otimes V^*$ are not isomorphic and some troubles appear with the 
 coproduct as  dual of the multiplication map. The concept of pairing solves these
difficulties.

A pairing between two Hopf algebras \cite{ChP},
$H$ and $H'$, is a bilinear mapping
$\langle \, \cdot \, ,  \, \cdot \, \rangle:
H \times H' \rightarrow \K$ that verifies the following
properties:
\begin{equation}\label{pairing}
\begin{array}{c}
\begin{array}{cc}
\langle h , m'(h'\otimes k')\rangle=
\langle \Delta(h) , h' \otimes k' \rangle, &  \langle h, 1_{H'} \rangle=
\epsilon(h),\\[2mm]
\langle h\otimes k , \Delta'(h')\rangle=
\langle m(h\otimes k) , h' \rangle, &  \epsilon'(h')= \langle  1_{H}, h'
\rangle,\end{array}\\ \\[-3mm]
\langle h, S'(h') \rangle =     \langle S(h), h' \rangle.
\end{array}
\end{equation}
Remark that  $
\langle h\otimes k, h'\otimes k' \rangle=
\langle h, h' \rangle  \langle k, k' \rangle.$

The pairing is said to be left (right) nondegenerate  if
$[\langle h, h' \rangle = 0, \quad \forall h' \in H'] \Rightarrow h=0$
($[\langle h, h' \rangle = 0, \quad \forall h \in H] \Rightarrow h'=0$). If the  pairing
 is simultaneously  left and right nondegenerate we simply say that it is
nondegenerate. 

The triplet  $(H, H', \langle\, \cdot \, , \, \cdot \, \rangle)$  composed by
two Hopf algebras and a  nondegenerate  pairing  will be called a `nondegenerate
triplet'. 

The bases $(h^m)$ of $H$ and
$(h'_n)$ of $H'$ are said to be dual with respect to the nondegenerate pairing if
\begin{equation}
\langle h^m, h'_n \rangle = c_n \delta^{m}_{n},\qquad  c_n \in \K -\{ 0 \}.
\end{equation}
The map   $f^\dagger: {H'} \rightarrow { H'}$ implicitly defined  in terms of 
the  map $f: H \rightarrow H$ by
\begin{equation}
\langle h, f^\dagger(h') \rangle= \langle f(h), h' \rangle ,
\end{equation}
is called the adjoint map of $f$ with respect to the nondegenerate 
pairing.

\subsect{Modules and comodules}
\label{modules}

Let us consider the triad $(V, \alpha, A)$,  where
$A$ is an associative  $\mathbb{K}$--algebra
with unit, $V$ is a  $\mathbb{K}$--vector space and  $\alpha$ a
linear map, $\alpha: A\otimes_{\mathbb{K}} V \rightarrow V$, called  action
and denoted  by  $a \lact v =\alpha(a\otimes v)$. We will say that  $(V, \alpha, A)$ 
(or $(V, \lact, A)$) is a left $A$--module if the following two conditions are verified:
\begin{equation}
    a \lact (b\lact v)= (ab) \lact v, \qquad 1 \lact v= v, \qquad
         \forall a, b \in A, \quad \forall v \in V.
\end{equation}

A morphism of left  $A$--modules, $(V, \lact,A)$ and $(V', \lact',A)$, is a linear
map, $f:V\rightarrow V'$,  equivariant with respect the action, i.e., 
\begin{equation} 
f(a\lact v)= a \lact' f(v), \qquad \forall a \in A,\ \forall v \in V. 
\end{equation}

Dualizing the concept of $A$--module it is obtained the concept of comodule.
Thus, if $C$ is an associative $\mathbb{K}$--coalgebra  with counit,
$V$ a  $\mathbb{K}$--vector space and 
$\beta: V \rightarrow C \otimes_{\mathbb{K}} V$ a linear map that will be called 
coaction and denoted by $v \lcact= \beta(v)  = v^{(1)} \otimes v^{(2)}$,  the triad  
$(V, \beta, C)$ (or $(V, \lcact, C)$) is said to be a 
 left $C$--comodule if the following axioms are verified:
  \begin{equation} \label{comodulo}
 {v^{(1)}}_{(1)} \otimes {v^{(1)}}_{(2)} \otimes {v^{(2)}}= 
           {v^{(1)}} \otimes {v^{(2)}}_{(1)} \otimes {v^{(2)}}_{(2)}, \quad
           \epsilon(v^{(1)}) v^{(2)}= v, \qquad \forall v \in V,
\end{equation}
where the coproduct of the elements of $C$ is symbolically written as 
$\Delta(c)=c_{(1)} \otimes c_{(2)}$.

  A linear map $f:V\rightarrow V'$ between
two  $C$--comodules, $(V, \lcact ,C)$ and  
$(V',\lcact',C)$ is a morphism if
\begin{equation} \label{comodulomorfismo}
   v^{(1)} \otimes f(v^{(2)})=f(v)^{(1)'}\otimes  f(v)^{(2)'},
   \qquad  \forall v \in V.
\end{equation}

Similarly right  $A$--modules and right $C$--comodules are defined.

\subsect{Module Algebras}

When a bialgebra acts or coacts on a vector space equipped with an additional  structure
of algebra, coalgebra or bialgebra \cite{Molnar,Maj95a} it is usual to demand some 
compatibility relations for the action. In the following  $B$ and  $B'$ will denote
bialgebras,  $A$ an algebra and $C$ a coalgebra.

The left module  $(A, \lact, B)$ is said to be a $B$--module algebra if
$m_A$ and $\eta_A$ are morphisms of $B$--modules, i.e.,  if 
\begin{equation} 
b\lact( aa')= (b_{(1)}\lact a)(b_{(2)}\lact a'), \qquad 
b \lact 1 = \epsilon(b) 1, \qquad 
\forall b \in B, \, \forall a, a' \in A.
\end{equation}

Changing algebra by  coalgebra it is obtained the structure of module coalgebra.
In this case the left $B$--module  $(C, \lact, B)$ is a $B$--module coalgebra if
$\Delta_C$ and $\epsilon_C$ are morphisms of $B$--modules, i.e., if 
$$
 (b\lact c)_{(1)} \otimes (b \lact c)_{(2)} =
    (b_{(1)}\lact c_{(1)}) \otimes (b_{(2)}\lact c_{(2)}), 
  \quad \epsilon_{C}(b \lact c) = \epsilon_{B}(b) \epsilon_{C}(c),
\qquad \forall b,c \in B.
$$

Dualizing these  structures two new ones are obtained.
The left $B$--comodule  $(C, \lcact, B)$ is said to be a  $B$--comodule coalgebra if
$\Delta_C$ and $\epsilon_C$ are morphisms of $B$--comodules, i.e., 
$$
    c^{(1)} \otimes {c^{(2)}}_{(1)} \otimes {c^{(2)}}_{(2)}=
   {c_{(1)}}^{(1)} {c_{(2)}}^{(1)} \otimes {c_{(1)}}^{(2)}
    \otimes {c_{(2)}}^{(2)}, \quad
     c^{(1)} \epsilon_C (c^{(2)})= (\eta_{B} \circ \epsilon_C)(c).
$$

The left $B$--comodule  $(A, \lcact, B)$ is a $B$--comodule algebra if
$m_A$ and $\eta_A$ are morphisms of $B$--comodules. Explicitly 
\begin{equation}
(aa')_{(1)} \otimes (aa')_{(2)} = a_{(1)}a'_{(1)} \otimes  a_{(2)} a'_{(2)}, \quad  
1_A \lcact = 1_B \otimes 1_A.
\end{equation}

The triad $(B', \lact, B)$ is a left $B$--module  bialgebra if simultaneously is a
$B$--module algebra and a $B$--module coalgebra; $(B', \lcact, B)$ is a left
$B$--comodule bialgebra if simultaneously is a $B$--comodule algebra and a
$B$--comodule coalgebra.

The  corresponding versions at the right are defined in an analogous manner.

By regular module (comodule) we understand an $A$--module ($C$--comodule) whose
vector space is the underlying vector space of the algebra $A$ (coalgebra $C$). The
action (coaction) is defined by means of the algebra product (coalgebra coproduct).  

For instance,   on the regular $A$--modules $(A, \lact, A)$ and $(A, \ract, A)$ the
actions are, respectively,
\begin{equation}
                 a \lact a'= aa', \qquad a' \ract a= a'a,
\end{equation}

If $B$ is a bialgebra,  the regular $B$--module $(B, \lact, B)$ whose `regular' action is defined
by
 \begin{equation} b \lact b' =b b',
\end{equation}
 is a module coalgebra. The module  $(B^*, \ract, B)$, obtained by dualization,
 is a module algebra with the `regular' action   
\begin{equation} 
\varphi \ract b = \langle \varphi_{(1)},b \rangle \varphi_{(2)}, \qquad b 
\in B, \; \varphi \in B^*. 
\end{equation}
 It will be also called regular module. The comodule versions can be easily obtained by
the reader.
  
\subsect{Bicrossproduct structure}
\label{bicross}

The concepts of  module algebra and  comodule coalgebra allow to describe
in a suitable way `semidirect' structures \cite{Molnar,Maj95a} as  we shall see
later.

Let $H$ be a bialgebra and  $(A, \ract , H)$  a right $A$--module  algebra. The 
expression
 \begin{equation}
    (h \otimes a) (h' \otimes a')= h h'_{(1)} \otimes (a \ract h'_{(2)}) a'
 \end{equation}
defines an algebra structure over $H \otimes A$,  denoted by
$H \RIMO A$ and called  semidirect product at the right (or simply right semidirect
product) of $A$ and $H$.

The `left' version is as follows:
let $(A, \lact , H)$ be a left $A$--module algebra.
A structure of algebra over $H \otimes A$, denoted by 
$A \LEMO H$ and called left semidirect product  
of $A$ and $H$, is defined by means of 
 \begin{equation}
(a \otimes h) (a' \otimes h')= a (h_{(1)} \lact a') \otimes h_{(2)} h' .
 \end{equation}

Dual structures  of the above ones are constructed in the following way.
Let $(C, \lcact, H)$ be a left $C$--comodule coalgebra. A coalgebra structure over 
$C\otimes H$, denoted by $C
\leco H$ and called  left semidirect product, is obtained  if
 \begin{equation} \begin{array}{lll}   
\Delta(c \otimes h) &=& c_{(1)} \otimes {c_{(2)}}^{(1)} h_{(1)}
             \otimes {c_{(2)}}^{(2)} \otimes h_{(2)}, \\[0.2cm]
        \epsilon(c \otimes h)&=& \epsilon_C(c) \epsilon_H(h).
\end{array}\end{equation}

 When $(C, \rcact, H)$ is a right $C$--comodule coalgebra,
the expressions 
 \begin{equation}\begin{array}{lll}     
\Delta(h \otimes c) &=& h_{(1)} \otimes {c_{(1)}}^{(1)} \otimes
h_{(2)}{c_{(1)}}^{(2)} \otimes c_{(2)},  \\[0.2cm]
\epsilon(h \otimes c) &=& \epsilon_C(h) \epsilon_H(c) , 
\end{array}\end{equation}
characterize a coalgebra structure over  $C \otimes H$ denoted by
 $ C \RICO H $ and called right semidirect product 
of $C$ and $H$. 

Let $K$ and $L$ be two bialgebras, such that  
$(L, \ract, K)$ is  a right  $K$--module algebra and $(K,\lcact,L)$ a left
$L$--comodule coalgebra. The tensor product
$K\otimes L$ is equipped simultaneously with the semidirect  structures of  algebra
$K\RIMO L$ and coalgebra $K \leco L$. If the following compatible conditions
are verified 
 \begin{equation}\label{compatibilidadbicross}
\begin{array}{c}
\epsilon(l \triangleleft k)= \epsilon(l) \epsilon(k),  \quad
\Delta(l\triangleleft k)= (l_{(1)} \triangleleft k_{(1)})
{k_{(2)}}^{(1)}  \otimes l_{(2)} \triangleleft {k_{(2)}}^{(2)},
\\ [0.2cm]
1 \lcact = 1 \otimes 1, \quad  (kk') \lcact= (k^{(1)}\triangleleft
k'_{(1)}){k'_{(2)}}^{(1)} \otimes
k^{(2)}{k_{(2)}'}^{(2)},\\ [0.2cm] {k_{(1)}}^{(1)} (l
\triangleleft k_{(2)}) \otimes {k_{(1)}}^{(2)}=  (l \triangleleft k_{(1)})
{k_{(2)}}^{(1)} \otimes {k_{(2)}}^{(2)},
\end{array}
 \end{equation}
then $K \RIMO L$ and  $K \leco L$ determine a bialgebra called (right--left) 
bicrossproduct and denoted by  $K \RL L$. 

If $K$ and $L$ are two Hopf algebras then $K \RL L$ has also an antipode given by
\begin{equation} 
S(k \otimes l) = (1 \otimes S(k^{(1)}l)) (S(k^{(2)}) \otimes 1).
  \end{equation}

On the other hand, let $\cal K$ and $\cal L$ be two bialgebras and  
$({\cal L}, \lact,{\cal K})$ and $({\cal K}, \rcact,{\cal L})$
a left $K$--module algebra and a
 right $L$--comodule coalgebra, respectively, verifying
the compatibility  conditions
\begin{equation}
\begin{array}{c}
\v(\lambda\act \k)=\v(\lambda)\v(\k),\\ [0.2cm]
\Delta(\lambda\act \k)\equiv (\lambda\act \k)_{(1)}\otimes (\lambda\act \k)_{(2)}
=({\lambda_{(1)}}^{(1)}\act \k_{(1)})\otimes {\lambda_{(1)}}^{(2)}(\lambda_{(2)}
\act \k_{(2)}),\\ [0.2cm]
\rcact (1)= 1\otimes 1,\\ [0.2cm]
\rcact (\k\k')
={\k_{(1)}}^{(1)}\k^{'(1)}\otimes {k_{(1)}}^{(2)}(\k_{(2)}
\act \k^{'(2)}) ,\\ [0.2cm]
{\lambda_{(2)}}^{(1)}\otimes(\lambda_{(1)}\act \k){\lambda_{(2)}}^{(2)}=
{\lambda_{(1)}}^{(1)}\otimes {\lambda_{(1)}}^{(2)}(\lambda_{(2)}\act \k).
\end{array}
\end{equation}
Then ${\cal L} \LEMO {\cal K}$ and  ${\cal L} \RICO {\cal K}$ determine
a bialgebra called (left--right) bicrossproduct  denoted by
${\cal L} \LR {\cal K}$.  

If $\cal K$ and $\cal L$ are two  Hopf algebras then ${\cal L} \LR {\cal K}$
has an antipode defined by
  \begin{equation}
S(\lambda  \otimes \kappa) = (1 \otimes S \kappa^{(1)}) 
(S(\lambda \kappa^{(2)}) \otimes 1).  
\end{equation}

Note that both bicrossproduct structures are
related by duality. Effectively, it can be proved that if 
$K$ and $L$ are two finite dimensional bialgebras,  and
the right  $K$--module algebra $(L, \ract, K)$ and the left
$L$--comodule coalgebra $(K,\lcact,L)$  verify the 
conditions (\ref{compatibilidadbicross}), then  $(K\RL L)^*= K^* \LR L^*$.

\subsect{ Star structures over bicrossproduct Hopf algebras}
\label{starstructures}

The following  original results show how construct dual bases and $*$--structures over
Hopf algebras with the structure of bicrossproduct when the  corresponding objects for the 
factors of the bicrossproduct are known \cite{oscarTesis}.

\begin{theorem}
\label{tbasesduales}
Let $H=K \RL L$ be a Hopf algebra  with  structure of bicrossproduct, and
$\langle \cdot, \cdot \rangle_1$ and $\langle \cdot, \cdot \rangle_2$ nondegenerate
pairings for the pairs $(K,K^*)$ and $(L, L^*)$, respectively. Then the expression
 \begin{equation} \label{pairingprod}
     \langle kl, \kappa \lambda \rangle=
    \langle k, \kappa \rangle_{1}
\langle l, \lambda \rangle_{2}.
 \end{equation}
defines a nondegenerate pairing between $H$ and $H^*$.
\end{theorem}

\begin{proof}[Proof]
Firstly note that
\begin{equation}
  \langle 1, \kappa \lambda \rangle=
  \langle 1_K \otimes 1_L, \kappa \lambda \rangle =
  \langle 1_K,\kappa   \rangle \langle 1_L, \lambda \rangle =
  \epsilon(\kappa) \epsilon(\lambda)=
  (\epsilon\otimes \epsilon) (\kappa \otimes \lambda)=
\epsilon(\kappa \lambda).
\end{equation}
On the other hand
\begin{equation}
\begin{array}{rl}
\langle kl, (\kappa \lambda)(\kappa' \lambda') \rangle = &
\langle kl, \kappa (\lambda_{(1)} \lact \kappa') \lambda_{(2)} \lambda'\rangle =
\langle k,\kappa (\lambda_{(1)} \lact \kappa') \rangle_1
\langle l, \lambda_{(2)} \lambda' \rangle_2\\[2mm]
 =& \langle k_{(1)}, \kappa \rangle_1 \langle k_{(2)}, \lambda_{(1)} \lact \kappa'
\rangle_1
    \langle l_{(1)}, \lambda_{(2)} \rangle_{2} \langle l_{(2)}, \lambda' \rangle_2
\\[2mm]
 =& \langle k_{(1)}, \kappa \rangle_1
   \langle \tau(k_{(2)} \lcact), \kappa' \otimes  \lambda_{(1)} \rangle
    \langle l_{(1)}, \lambda_{(2)} \rangle_{2} \langle l_{(2)}, \lambda' \rangle_2
\\[2mm]
 =& \langle k_{(1)}, \kappa \rangle_1
    \langle {k_{(2)}}^{(2)}, \kappa' \rangle_{1}
   \langle {k_{(2)}}^{(1)},  \lambda_{(1)} \rangle_2
    \langle l_{(1)}, \lambda_{(2)} \rangle_{2}
     \langle l_{(2)}, \lambda' \rangle_2 \\[2mm]
 =& \langle k_{(1)}, \kappa \rangle_1
   \langle {k_{(2)}}^{(2)}, \kappa' \rangle_{1}
   \langle {k_{(2)}}^{(1)}l_{(1)},  \lambda \rangle_2
        \langle l_{(2)}, \lambda' \rangle_2 \\[2mm]
 =& \langle k_{(1)}{k_{(2)}}^{(1)}l_{(1)}, \kappa \lambda\rangle
    \langle {k_{(2)}}^{(2)} l_{(2)}, \kappa'\lambda' \rangle \\[2mm]
 =&\langle \Delta(kl), (\kappa \lambda) \otimes (\kappa' \lambda') \rangle.
\end{array}
\end{equation} 
Similarly for the identities
\begin{equation}
  \langle kl, 1 \rangle= \epsilon(kl), \qquad
  \langle (kl) (k'l'), \kappa \lambda \rangle =
  \langle (kl) \otimes (k'l'), \Delta(\kappa \lambda) \rangle.
\end{equation}
Hence, it is proved that
 $\langle \cdot, \cdot \rangle$ is a bialgebra pairing.
The pairing is nondegenerate. Effectively, fixed a basis 
$(l_i)_{i\in I}$ of $L$,
the coproduct can be written as
\begin{equation}
 \Delta(h)= \sum_{i\in I} a_i(h) \otimes l_i,
\end{equation}
with $a_i: H \to K$. 
Let us suppose that
\begin{equation} \label{productonulo}
    \langle h, \eta \rangle =0, \qquad  h\in H,\  \eta \in H^*.
\end{equation}
Expression (\ref{productonulo}) can be rewritten as
$$ 
\langle h, \kappa \lambda \rangle =\langle \Delta(h), \kappa\otimes
\lambda\rangle 
  =   {\displaystyle \sum_{i \in I}} \langle a_i(h),
    \kappa \rangle_1  \langle l_i , \lambda \rangle_2
   =   \langle  { \displaystyle \sum_{i \in I}} \langle  a_i(h),\kappa \rangle_1  
l_i ,  \lambda \rangle_2 =0,
$$
where $\kappa\in K^*, \ \lambda \in L^*$. Since $\langle \cdot, \cdot
\rangle_2$ is nondegenerate and
$\langle \cdot, \cdot \rangle_1$ is also nondegenerate  we get that
$a_i(h)=0,$ hence $\Delta(h)= 0$.
 Finally, using the counit axiom
\begin{equation}
    h= (\epsilon \otimes {\rm id}) \circ \Delta(h)= 0, 
\end{equation}
we have proved that the pairing is left nondegenerate. In a similar way it is proved  that the
pairing is nondegenerate at the right.  Using the fact that the last equality
of (\ref{pairing}) is a consequence of the two first ones when the pairing is
nondegenerate, we conclude that the bilinear form (\ref{pairingprod}), which is a
bialgebra pairing, is also a pairing of Hopf algebras.
\end{proof}
\begin{corollary}
  \label{basesdualesbicross}
With the  pairing and the notation defined in the previous theorem if   $(k_m)$ and
$(\kappa_m)$ are dual bases for $K$ and $K^*$, and $(l_n)$ and $(\lambda_n)$ are
dual  bases for $L$ and $L^*$, then
$(k_m l_n)$ and $(\kappa^m \lambda^n)$ are dual bases for $H$ and $H^*$. In other
words, if  
$  \langle k_m, \kappa^{m'} \rangle  = \delta_m^{m'} $
and $  \langle l_n, \lambda^{n'} \rangle = \delta_n^{n'}$
then    $ \langle k_m l_n, \kappa^{m'} \lambda^{n'} \rangle=
\delta_m^{m'}\delta_n^{n'}$.
\end{corollary}

In  the case  of left-right bicrossproduct there is  a similar result.

\begin{theorem}
\label{estructuraestrellabicross}
Let us consider the bicrossproduct Hopf algebra $H=K\RL L$. Supposing that
$K$ and $L$ are equipped with $*$--structures with the
following  compatibility relation
 \begin{equation} \label{compatibilidadestrellabicross}
         (l\ract k)^*= l^* \ract S(k)^*.
 \end{equation}
Then the expression
 \begin{equation} \label{estrellabicross}
     (kl)^*= l^* k^*, \qquad k \in K, \ l \in L,
  \end{equation}
determines a $*$--structure on the algebra sector of $H$.
\end{theorem}

\begin{proof}[Proof]
The definition of a $*$--structure on  $H$ has to be consistent
with  the algebra  structure  is an antimorphism, i.e.,
  \begin{equation} \label{inicio}
      (lk)^*= k^* l^*, \qquad k \in K, \quad l \in L.
  \end{equation}
Since the product  on $H$ establishes that
 \begin{equation}
        lk= k_{(1)} (l \ract k_{(2)}),
 \end{equation}
and according to the definition (\ref{estrellabicross}) 
 \begin{equation}
        (lk)^*= (l\ract k_{(2)})^* (k_{(1)})^*.
 \end{equation}
Using the product on $H$ one obtains
 \begin{equation} \label{mari}
(lk)^*= [(k_{(1)})^*]_{(1)} \; \{ (l\ract k_{(2)})^* \ract [(k_{(1)})^*]_{(2)}\}.
\end{equation}
Taking into account (\ref{compatibilidadestrellabicross}) and that the $*$--structure
on $K$ is a coalgebra morphism the  equality (\ref{mari}) becomes
 \begin{equation}
   (lk)^*= (k_{(1)})^*  \; \{l^* \ract [S(k_{(3)})^* k_{(2)}^*] \}.
 \end{equation}
Finally, the property characterizing the antipode reduces  this expression to 
(\ref{inicio}).
\end{proof}

\sect{Induced representations for quantum  algebras}
\label{induced representations}

Since the algebras involved in this work are equipped with a bicrossproduct
structure different actions appear. In order to avoid any confusion we will
denote them by the following symbols (or their symmetric for the corresponding
right actions and  coactions):

 \hspace{1cm} $\lact$  ($\lcact$): actions  (coactions)
      of the bicrossproduct structure,

 \hspace{1cm} $\vdash$: induced and inducting representations,

 \hspace{1cm} $\succ$ ($\prec$): regular actions (coactions).

In the following we will show that the problem of the determination of the induced
representations is reduced as a last resort to the
expression of  products in normal  ordering. The next result will be very useful for
this purpose.

\begin{proposition} \label{lemaconmutacion}
Let $A$ be an associative algebra. The following relations hold:
\begin{equation}\begin{array}{lll}\label{relacionesadjuntas}
        a^m a' &=& \sum_{k=0}^m \comb{m}{k} \lad{a}^k(a') a^{m-k}, \\[0.6cm]
        a' a^m &=& \sum_{k=0}^m \comb{m}{k}  a^{m-k} \rad{a}^k(a'), 
\end{array}\qquad \forall a,a' \in A,\ m \in  \mathbb{N},
     \end{equation}
\end{proposition}
where 
\begin{equation} \label{accionadjunta}
      \lad{a}(a')= aa' -a'a= [a,a'],\qquad
       \rad{a}(a')= a'a -aa'=[a',a].
\end{equation}
\begin{proof}[Proof]
The demonstration is by  induction. The relations
(\ref{relacionesadjuntas}) are trivial  identities for $m=0$. Let us suppose that
the first expression is true for $m\in \mathbb{N}$, then for $m+1$ we have
 \begin{equation} \begin{split}
    a^{m+1} a'=& a (a^m a') = a \sum_{k=0}^m \comb{m}{k} \lad{a}^k(a') a^{m-k} \\
    = & \sum_{k=0}^m \comb{m}{k} [\lad{a}^k(a') a + \lad{a}^{k+1}(a')] a^{m-k} \\
    = & \sum_{k=0}^m \comb{m}{k} \lad{a}^k(a') a^{m-k+1} +
          \sum_{k=0}^m \comb{m}{k} \lad{a}^{k+1}(a') a^{m-k} \\
    = & \sum_{k=0}^m \comb{m}{k} \lad{a}^k(a') a^{m-k+1} +
          \sum_{k=1}^{m+1}\comb{m}{k-1} \lad{a}^{k}(a') a^{m-k+1} \\
    = & \sum_{k=0}^{m+1} \comb{m+1}{k} \lad{a}^k(a') a^{m+1-k}.
    \end{split}
  \end{equation}
The proof of the second identity (\ref{relacionesadjuntas}b) is similar.
\end{proof}

Note that in an appropriate  topological context, where it is allowed the
convergence and the reordering of  series, expressions (\ref{relacionesadjuntas})
carry to the usual relation between adjoint action and
exponential mapping:
\begin{equation} \begin{split}
  e^a a' =& \sum_{m=0}^\infty \frac{1}{m!} \sum_{k=0}^m \comb{m}{k} 
\lad{a}^k(a') a^{m-k} = \sum_{m=0}^\infty \sum_{k=0}^m \frac{1}{k! (m-k)!}
\lad{a}^k(a') a^{m-k} \\
   = &    \sum_{k=0}^\infty \sum_{m=k}^\infty \frac{1}{k! (m-k)!} 
\lad{a}^k(a') a^{m-k}  =  e^{\lad{a}}(a')   e^{a},
 \end{split}
\end{equation}
or equivalently
\begin{equation} \label{lemconmexp}
  e^a a' e^{-a} = e^{\lad{a}}(a').
\end{equation}
For the other adjoint action  taking into account that $\lad{-a} =\rad{a}$ we get an
analogous relation 
\begin{equation}  \label{lemconmexp2}
      e^{-a} a' e^{a} =  e^{\rad{a}}(a').
\end{equation}
\subsect{General case}
\label{indgen}

Let us consider a nondegenerate triplet $(H, {\cal H}, \langle \cdot, \cdot\rangle)$.
Let $L$ be a commutative subalgebra of $H$ and $\{ l_1, \ldots,l_s\}$, a system of
generators of  $L$ which can be completed with
$\{ k_1,\ldots,k_r\}$, in such a way that  $(l_n)_{n \in \mathbb{N}^s}$ is a basis of
$L$ and $(k_ml_n)_{(m,n) \in \mathbb{N}^r \times \mathbb{N}^s}$ a basis of 
$H$. Moreover, suppose that there is a system of generators in $\cal H$,
$\{ \kappa_1, \ldots, \kappa_r, \lambda_1, \ldots, \lambda_s\}$, such that 
$(\kappa^m \lambda^n)_{(m,n) \in \mathbb{N}^r \times \mathbb{N}^s}$ is a basis
of $\cal H$ dual of that of  $H$ with the pairing
\begin{equation} 
\langle k_m l_n, \kappa^{m'} \lambda^{n'} \rangle = m! n! \;
\delta^{m'}_m \delta^{n'}_n.
\end{equation}

We are interested in the description of the representation induced by  the character
of $L$ determined by  $a=(a_1,\ldots, a_s) \in \mathbb{K}^s$, i.e., 
\begin{equation}
1 \dashv l_n = a_n=
a_1^{n_1} \cdots a_s^{n_s}, \qquad n \in  \mathbb{N}^s.
\end{equation}
The  elements $f$ of $\text{Hom}_\mathbb{K}(H,\mathbb{K})$ verifying the
invariance  condition
\begin{equation}\label{equivariancecondition}
   f(hl)= f (h) \dashv l, \qquad
    \forall l \in L, \quad \forall h \in H,
\end{equation}
constitute the carrier space $\mathbb{K}^\uparrow =\text{Hom}_L(H,\mathbb{K})$ of
the induced representation. 
Identifying $\text{Hom}_\mathbb{K}(H,\mathbb{K})$ with
$\cal H$ using the pairing, the elements of
$f \in \mathbb{K}^\uparrow$  can be written as
\begin{equation}
f= \sum_{(m,n)\in \mathbb{N}^r \times \mathbb{N}^s} f_{m n}
 \kappa^{m} \lambda^{n}.
\end{equation}
The equivariance condition (\ref{equivariancecondition})
\begin{equation}
   \langle hl, f\rangle = \langle h, f \rangle \dashv l, \qquad
    \forall l \in L, \quad \forall h \in H,
\end{equation}
combined with duality  gives the following relation between the 
coefficients $f_{m n}$
\begin{equation}
     m! n! f_{m n}= \langle k_m l_n, f \rangle =
     \langle k_m, f \rangle  a_n = m! f_{m 0} a_n.
\end{equation}
Hence, the elements of the carrier space of the induced representation are
\begin{equation}
f= \kappa \psi, \qquad \kappa \in {\cal K},
\end{equation}
where $\psi=e^{a_1 \lambda_1} \cdots e^{a_s \lambda_s}$,  and
$\cal K$ is the  subspace of $\cal H$ generated by the linear  combinations
of the ordered monomials $(\kappa^m)_{m \in \mathbb{N}^r}$. Since 
$\psi$ is invertible (it is product of exponentials) there is an isomorphism between
the vector spaces ${\cal K}$ and $\mathbb{K}^\uparrow$ given by 
$\kappa  \to  \kappa \psi$.

The action of $h \in H$ over the elements of  $\mathbb{K}^\uparrow$ is determined
knowing the action over the basis elements $(\kappa^p \psi)_{p \in \mathbb{N}^r}$ of
this space. So, putting
\begin{equation}
(\kappa^p \psi) \dashv h =
   \sum_{(m,n)\in \mathbb{N}^r \times \mathbb{N}^s}  [h]^p_{m n} \kappa^m\lambda^n, 
\qquad  p \in \mathbb{N}^r,
\end{equation}
the constants $[h]^p_{m n}$ can be evaluated by means of duality
 \begin{equation}\label{coef} 
m! n! [h]^p_{m n} = \langle (\kappa^p \psi) \dashv h,
 k_m l_n\rangle =
    \langle \kappa^p \psi,  h k_m l_n\rangle =
    \langle \kappa^p \psi,  h k_m \rangle a_n.
 \end{equation}

The properties of the action allow to compute it only for the generators of $H$ instead of considering
an arbitrary element $h$ of $H$. Finally, all that  reduces to write the product $h k_m$ in normal
ordering to get the value of the paring in (\ref{coef}). However in many cases this task is very
cumbersome, for this  reason now our objective is to take advantage of the bicrossproduct structure to
simplify the job.

\subsect{Quantum algebras with bicrossproduct structure}
\label{inducedrepresentationsbicross}

In the following we will restrict ourselves to Hopf algebras having a
bicrossproduct structure like  $H= {\cal K} \bicross {\cal L}$, such 
that the first factor is cocommutative and the second commutative.

We are interested in the construction of the representations induced by 
`real' characters of the commutative sector  $\cal L$.
We will show that the solution of this problem can be reduced to the study of
certain dynamical systems which  present, in general, a non linear action.

Let us start adapting the construction presented in the previous subsection
\ref{indgen} to the  bicrossproduct  Hopf algebras
$H= {\cal K}\bicross {\cal L}$.
Let us suppose that the  algebras $\cal K$ and $\cal L$ are finite 
generated by the sets $\{ k_i \}_{i=1}^r$ and  
$\{ l_i \}_{i=1}^s$, respectively, such that the generators $k_i$
are primitive.  

We also assume that  $(k_n)_{n\in \mathbb{N}^r}$ and $(l_m)_{m\in \mathbb{N}^s}$
are bases of the vector spaces underlying   to $\cal K$
and $\cal L$, respectively.
Let ${\cal K}^*$ and ${\cal L}^*$ be the dual algebras of $\cal K$ and $\cal L$,
respectively,   having  dual systems to those of $\cal K$ and $\cal L$  with analogue
properties  to them. Hence,  duality between  $H$ and $H^*$ is given by  
\begin{equation}\label{dualidad}
     \langle k_m l_n, \kappa^{m'} \lambda^{n'} \rangle =
     m! n! \; \delta_m^{m'} \delta_n^{n'}.
\end{equation}
As we will see later these hypotheses are not, in reality,  too restrictive.
All these generator systems will be used to described the induced
representations. 

Let us consider the  character  of  $\cal L$ labeled  by
$a \in \mathbb{C}^s$ 
\begin{equation}  \label{caracter}
 1 \dashv l_n= a_n, \qquad n \in \mathbb{N}^s,
\end{equation}
the  discussion of subsection \ref{indgen} allows us to state the
following theorem.
\begin{theorem}\label{tind0} 
The   carrier space, $\mathbb{C}^\uparrow$,
of  the  representation of $H$ induced by the  character $a$ of  $\cal L$ (see
(\ref{caracter})) is isomorphic to ${\cal K}^*$ and is constituted by the  elements of  the  form  
$\kappa \psi$, where
$\kappa \in {\cal K}^*$ and 
\begin{equation}
         \psi= e^{a_1 \lambda_1} e^{a_2 \lambda_2} \cdots
              e^{a_s \lambda_s}.
\end{equation}
The  induced action  is given by 
\begin{equation} \label{accion0}
           f \dashv h= \sum_{m\in \mathbb{N}^r}
            \kappa^m \langle h \frac{k_m}{m!}, f \rangle \psi,
            \qquad h\in H,\  f\in \mathbb{C}^\uparrow.
\end{equation}
\end{theorem}

The  action  of the generators of $\cal K$ and $\cal L$ in  the induced 
representation   will be given in the next theorem, which needs the introduction of 
some new concepts.

Since $\cal L$ is commutative, it can be identified with the algebra of
functions $F(\mathbb{R}^s)$ by means of the  algebra morphism
${\cal L} \to F(\mathbb{R}^s),\ (l \mapsto \tilde{l})$,
which maps the  generators of  $\cal L$ into the canonical projections
\begin{equation}
 \tilde{l}_j(x)= x_j, \qquad\quad 
\forall x=  (x_1, x_2, \ldots, x_s)\in \mathbb{R}^s, \quad
 j=1,2,\ldots, s.
\end{equation}
The  $*$--structure keeping
invariant the  generators chosen in $\cal L$ is distinguished  in a
natural way  by the above identification 
\begin{equation} \label{estre}
       l_j^*= l_j, \qquad j=1,2,\ldots, s.
\end{equation}
The characters (\ref{caracter}) compatible with (\ref{estre}) are `real', i.e., 
determined by the elements $a\in \mathbb{R}^n \subset \mathbb{C}^n$. We
will restrict to them henceforth.
Note that the  character    (\ref{caracter}), with $a \in \mathbb{R}^n$,
 can be written now as
\begin{equation}
  1 \dashv l = \tilde{l}(a).
\end{equation}
The  right action  of $\cal K$ on $\cal L$ can be translated to 
 $F(\mathbb{R}^s)$ because the generators of $\cal K$ are primitive and, hence, they
act by derivations on the   $\cal K$--module algebra of ${\cal K} \RL {\cal L}$.
Thus, the generators $k_i$ induce vector fields, $X_i$, on $\mathbb{R}^s$
determined by 
\begin{equation}
 X_i\, \tilde{l}= \widetilde{l\ract k_i},
   \qquad i=1,2, \ldots, r.
\end{equation}
The  corresponding  flow,
 $\Phi_i:\mathbb{R} \times \mathbb{R}^s \rightarrow \mathbb{R}^s$,
is implicitly defined by
\begin{equation} \label{flujocampo}
   (X_i f)(x)= (D f_{x,\Phi_i})(0),
\end{equation}
where $f_{x,\Phi_i}(t)= f\circ \Phi_i^t(x)$ and $D$ is the  derivative  operator
over real variable functions. Notice that, in 
general, the    one-parameter group of  transformations associated to the flow
$\Phi_i$ is not  globally defined. 

\begin{proposition}\label{lemaconmutacionbicross}
In the   Hopf algebra $H= {\cal K}\bicross {\cal L}$ the following relation holds
\begin{equation}
   l k_m= \sum_{p \leq m} \comb{m}{p} k_{m-p} (l \LLL k_p),
   \qquad  \forall l \in {\cal L}, \forall m \in \mathbb{N}^r,
\end{equation}
\noindent where the multi-combinatorial number is defined as product of usual
combinatorial numbers or through  multi-factorials
\begin{equation}
\comb{m}{p}= \prod_{i=1}^r \comb{m_i}{p_i} = \frac{m!}{p! (m-p)!},
\end{equation} 
where  the ordered  relation over the multi-indices is given by
\begin{equation}
p\leq m \, \, \, \Leftrightarrow \, \, \, p_1 \leq m_1,\  p_2 \leq m_2,
\ldots, \ p_r \leq m_r,
\end{equation}
and if $p \leq m$  the  difference  between $m$ and $p$ is well defined in
$\mathbb{N}^r$ by
\begin{equation}
 m-p=(m_1-p_1, m_2-p_2, \ldots, m_r-p_r).
\end{equation} 
\end{proposition}

\begin{proof}[Proof]
Let us consider an element $l$ of $\cal L$ and a generator
$k_i$ of  $\cal K$ in the associative  algebra ${\cal K}\RL {\cal L}$. Taking into
account  the  definition of the product in
${\cal K} \RL {\cal L}$ and that the generators
 $k_i$ are primitive we can write
\begin{equation}
     \rad{k_i}(l)= [l,k_i]= l \ract k_i, \qquad
     \left(\rad{k_i}\right)^p (l)= l \ract k_i^p.
\end{equation}
Picking out  the  second formula of 
(\ref{relacionesadjuntas}) for
$a'=l$ and $a=k_i$ we get
\begin{equation}
   l k_i^m= \sum_{p \leq m } \comb{m}{p} k_i^{m-p} (l \ract k_i^p).
\end{equation}
This formula is valid for $m\in \mathbb{N}$.
The  validity of the  expression for a  multi-index
 $m \in \mathbb{N}^r$ is a direct consequence of the  properties of  the  action 
$\ract$ and of the definitions of the  multi-objects that has been  introduced.
\end{proof}
\begin{theorem}\label{tind1}
The  explicit action  of the generators of  $\cal K$ and $\cal L$ in
 the  induced representation   of Theorem \ref{tind0}
realized in the   space ${\cal K}^*$ is given by the following  expressions:
\begin{equation} \label{accion1}
\begin{array}{rcl}
           \kappa \dashv k_i & = & \kappa \prec k_i \  ,\\[0.2cm]
           \kappa \dashv l_j & = & \kappa \, \,
        \hat{l}_j  \! \circ \Phi_{(\kappa_1,\kappa_2, \ldots, \kappa_r)}(a) \ ,
\end{array}
\end{equation}
where $i \in \{1,\ldots, r\}$, $j \in \{1,\ldots, s\}$, the   symbol
$\prec$ denotes the  regular action of  $\cal K$ on ${\cal K}^*$, and
$\Phi_{(\kappa_1,\kappa_2, \ldots, \kappa_r)}=
\Phi_r^{\kappa_r} \circ \cdots \circ   \Phi_2^{\kappa_2} \circ 
\Phi_1^{\kappa_1}$.
\end{theorem}

\begin{proof}[Proof]
For the  first expression we apply
(\ref{accion0}) to the case $h= k_r$
\begin{equation}\begin{array}{lll}
  (\kappa \psi) \dashv k_i&= & \sum_{m\in \mathbb{N}^r}
     \kappa^m \langle k_i \frac{k_m}{m!}, \kappa \psi \rangle \psi=
     \sum_{m\in \mathbb{N}^r}
   \kappa^m \langle k_i \frac{k_m}{m!}, \kappa  \rangle
           \langle 1_{\cal L}, \psi \rangle \psi \\[0.25cm]
  &=&  \sum_{m\in \mathbb{N}^r}
   \kappa^m \langle \frac{k_m}{m!}, \kappa \prec k_i \rangle
    \psi=(\kappa \prec k_i)   \psi.
  \end{array}
\end{equation}
For  the  third equality we use that 
$\langle 1_{\cal L}, \psi \rangle=1 $, and the  last one is based on the  fact
that $\frac{1}{m!}\kappa^m \otimes k_m$ is  the  $T$--matrix \cite{tmatrix} of the pair 
$({\cal K}^*, {\cal K})$. 

The   computation  of  the  action of  $l_j$ is more complicated
\begin{equation} \label{accion1lj}
\begin{array}{lll}
 (\kappa \psi) \dashv l_j&=&
    {\displaystyle \sum_{m \in \mathbb{N}^r}}
   \kappa^m \langle l_j \frac{\kappa_m}{m!}, \kappa \psi \rangle \psi
      \\[0.2cm]
 &= & {\displaystyle \sum_{m \in \mathbb{N}^r}}
    \kappa^m \langle {\displaystyle \sum_{p \leq m}} \comb{m}{p} \frac{1}
{m!} k_{m-p} (l_j  \ract k_p), \kappa \psi \rangle \psi    \\
 &= & {\displaystyle \sum_{m \in \mathbb{N}^r}} {\displaystyle \sum_{p \leq m}} 
\frac{1}{p! (m-p)!}
    \kappa^m \langle k_{m-p} (l_j  \ract k_p), \kappa \psi \rangle \psi
     \\[0.2cm]
 &= & {\displaystyle \sum_{p \in \mathbb{N}^r}} {\displaystyle \sum_{m \in p +
\mathbb{N}^r}}
    \frac{1}{p! (m-p)!}
    \kappa^m \langle k_{m-p} (l_j  \ract k_p), \kappa \psi \rangle \psi
      \\[0.2cm]
& = & {\displaystyle \sum_{p \in \mathbb{N}^r}} {\displaystyle \sum_{m \in 
\mathbb{N}^r}}
    \frac{1}{p! m!}
    \kappa^{m+p} \langle k_{m} (l_j  \ract k_p), \kappa \psi \rangle \psi
   \\[0.2cm]
  &= &   {\displaystyle \sum_{m \in  \mathbb{N}^r}} \frac{1}{m!} \kappa^{m}
          \langle k_{m}, \kappa \psi \rangle 
       {\displaystyle \sum_{p \in \mathbb{N}^r}}  \frac{1}{p!} \kappa^{p}
          (1 \dashv (l_j  \ract k_p)) \psi
    \\[0.2cm]
&  = &   {\displaystyle \sum_{m \in  \mathbb{N}^r}} \frac{1}{m!} \kappa^{m}
          \langle k_{m}, \kappa  \rangle 
       {\displaystyle \sum_{p \in \mathbb{N}^r}}  \frac{1}{p!} \kappa^{p}
          (1 \dashv (l_j  \ract k_p)) \psi
      \\[0.2cm]
&  = &   \kappa
       {\displaystyle \sum_{p \in \mathbb{N}^r}}  \frac{1}{p!} \kappa^{p}
          (1 \dashv (l_j  \ract k_p)) \, \psi
       \\[0.2cm]
 & = &   \kappa
       {\displaystyle \sum_{p \in \mathbb{N}^r}}  \frac{1}{p!} \kappa^{p}
          \, \widehat{l_j  \ract k_p}(a)  \, \psi
        \\[0.2cm]
 & = &   \kappa
       \left[ {\displaystyle \sum_{p \in \mathbb{N}^r}}  \frac{1}{p!} \kappa^{p}
           X_p\Bigg\vert_a \hat{l}_j  \right] \, \psi.
\end{array}
\end{equation}
In  the  second equality of (\ref{accion1lj}) Proposition
\ref{lemaconmutacionbicross} has been used. The next three are
simple reorderings of  the  sums.  The  sixth equality is a consequence of 
the  equivariance property  and of  the  commutativity in ${\cal K}^*$. The 
definitions of  the    duality form in  the  bicrossproduct structure, 
 of  the $T$--matrix  of the algebra $\cal K$,  of the 
identification of  $L$ with the   algebra of  functions $F(\mathbb{R}^s)$ are
 successively applied in the next equalities.  Finally, it is defined
$X_p= X_r^{p_r} \ldots X_2^{p_2} X_1^{p_1}$ in terms of the 
vector fields associated to the generators $k_i$.

On the other hand, from relation  (\ref{flujocampo}) between the
flow $\Phi_i$ and $X_i$ one gets
\begin{equation} \label{4f}
    f \circ \Phi_i^t(x)= f_{x,\Phi_i}(t)= (e^{tD} f_{x,\Phi_i})(0)=
     \sum_{n=0}^\infty \frac{1}{n!} t^n (D^n f_{x,\Phi_i})(0)=
      \sum_{n=0}^\infty \frac{1}{n!} t^n (X^n_i f)(x),
\end{equation}
for any regular function $f\in F(\mathbb{R}^s)$. So, to get the expression of the
action established in the theorem it suffices to take  
 $f=\hat{l}_j$,
$x=a$, and replacing formally the real number $t$ by $\kappa^i$, making
successively
 $i=1,\ldots r$, and substitute
 the  relation obtained in (\ref{accion1lj}).
\end{proof} 

Remark that the inverse order in  product $X_p= X_r^{p_r} \ldots
X_2^{p_2} X_1^{p_1}$ and in the flow composition
$\Phi_{t_1,t_2,\ldots,t_r}=\Phi_{t_r} \ldots \Phi_{t_2} \Phi_{t_1}$
is due to that the action of  $\cal K$ on $\cal L$ is at right.

When  $H$ is the deformed enveloping algebra  of  a semidirect  product
with  Abelian kernel and the   sector $\cal K$ is nondeformed then  the  first
expression  of  (\ref{accion1}) says that  the  representation of  $\cal K$ is
the  same that in the nondeformed case. On the other hand, the generators of 
$\cal L$ act as multiplication operators affected by the  deformation.

\sect{Modules and representations}
\label{modulesandrepresentations}

The deep relationship between representations and modules (see \cite{olmo00})  allow to
reformulate the theory of induced representations for quantum algebras that we have
developed in the previous section from the perspective of module theory.
\subsect{Regular modules}
\label{regularmodules} 

The   objective of  this section is to describe the four regular $H$--modules associated to a
Hopf algebra  $H$: $(H,\prec, H)$, $(H^*,\succ, H)$, $(H,\succ, H)$ and 
$(H^*,\prec,H)$; $H^*$ is the dual of $H$ in the sense of nondegenerate pairing (see subsection
\ref{duality}).  

It is well known the existence of theorems proving that,
essentially, all  the commutative or cocommutative Hopf algebras are of the form
$F(G)$ or $\mathbb{K}[G]$ (or $U(\mathfrak{g})$) for any  group $G$ \cite{Maj95a}.
So, the  kind  of  bicrossproduct that we will consider can be described as  
 \begin{equation}
     H= \mathbb{C}[K] \RL F(L) \qquad \text{or} \qquad
     H= U(\mathfrak{k})  \, \RL F(L),
 \end{equation}
where  $K$ and $L$ are finite groups or  Lie groups.

We will focus our attention in the case that both, $K$ and $L$,
are Lie groups with associated  Lie algebras $\mathfrak{k}$ and
$\mathfrak{l}$, respectively. In this way, the  dual of  $H$ will be
 \begin{equation}
    H^*= F(K) \LR  \, U(\mathfrak{l}).
 \end{equation}

The  clue for an effective  description  of  the regular
modules is the use of elements of  $H$ and $H^*$ like 
\begin{equation}\begin{array}{lllll}  \label{tablita}
&  k \lambda \in H , &  \qquad\qquad  & k \in K ,\ \  &\lambda \in F(L) ,\\[0.2cm]
& \kappa l \in H^* ,  &\qquad\qquad &\kappa \in F(K) , \ \  & l \in L .
\end{array}\end{equation}
We will see that these elements describe
 completely the structures of  the regular $H$--modules  and
are more  convenient than  the  bases of  ordered  monomials. 

\begin{theorem}
\label{crbicross0}
Let us consider  elements 
$k,k' \in K$, $\lambda,\lambda' \in F(L)$, $\kappa \in F(K)$ and $l \in L$.
The   action  on any of the four
regular $H$--modules is:
\begin{equation} \begin{array}{lllll}
(H,\prec, H)&:&\ \
 (k\lambda) \prec k'=  kk'(\lambda \ract k'), 
\qquad\quad  &(k\lambda) \prec \lambda' = k\lambda \lambda' \ ;\\[0.3cm]
(H^*,\succ, H)&:& \ \
 k' \succ (\kappa l)  =  (k'\succ \kappa) (k'\lact l),\qquad\quad  
&\lambda' \succ (\kappa l) = \kappa (\lambda' \succ l)\  ;\\[0.3cm]
(H,\succ, H)&:& \ \
 k' \succ (k\lambda) =  k'k\lambda ,\qquad\quad 
&\lambda' \succ (k\lambda) =  k (\lambda' \ract k) \lambda \ ;\\[0.3cm]
(H^*,\prec, H)&:& \ \
(\kappa l) \prec k'  =  (\kappa \prec k ') l , \qquad\quad
&(\kappa l) \prec \lambda'  = \kappa \langle l^{(1)}, \lambda' \rangle l^{(2)} l\ .
\end{array}\end{equation} 
\end{theorem}
\begin{proof}[Proof]
{\bf (1)} The results relative to  the modules
 $(H,\prec, H)$ and $(H,\succ, H)$
only require the use of  the  product
 defined on the   semidirect product of  algebras $U(\mathfrak{k})
\rimo F(L)$ (remember that for  arbitrary elements 
$k,k'\in U(\mathfrak{k})$ and $\lambda, \lambda' \in F(L)$ 
such product is given by 
$(k\otimes \lambda)(k'\otimes \lambda')=
k k'_{(1)} \otimes (\lambda \ract k'_{(2)}) \lambda'$).
In order to evaluate  the  action of  $k'$
we take into account that $\Delta(k')= k' \otimes k'$.

{\bf (2)} In the   module algebra  $(H^*,\succ, H)$
 the  action of  $k'$ is obtained by 
\begin{equation} \label{cadig1}
 \begin{array}{lll}
    \langle k \lambda, k' \succ (\kappa l)\rangle &= &
            \langle (k \lambda) \prec k', \kappa l \rangle \quad =
            \langle k k' (\lambda \ract k'), \kappa l \rangle  \\[0.2cm]
   & = &      \langle k k', \kappa \rangle \langle \lambda \ract k',  l \rangle\ =
            \langle k, k' \succ \kappa \rangle \langle \lambda,   k' \lact  l
\rangle \\[0.2cm]
   & = &       \langle k \lambda, (k' \succ \kappa) (k' \lact  l) \rangle.
   \end{array}
 \end{equation}
The  action of  $l'$ is obtained in an  analogous way
\begin{equation} \label{cadig2}
\begin{array}{lll}
  \langle k \lambda,  \lambda' \succ (\kappa l)\rangle &=&  
  \langle k \lambda \prec  \lambda', \kappa l\rangle =
    \langle k \lambda \lambda', \kappa l\rangle   \\[0.2cm]
 & = & \langle k, \kappa \rangle  \langle  \lambda \lambda',  l\rangle \ =
\langle k, \kappa \rangle  \langle  \lambda,  \lambda'\succ  l\rangle
  \\[0.2cm]
  &= & \langle k  \lambda,  \kappa (\lambda'\succ  l)\rangle.
   \end{array}
   \end{equation}
Notice that in the  first expression of  (\ref{cadig1}) and of
(\ref{cadig2}) the   symbol $\succ$ represents
 the  regular action  of  $(H^*,\succ, H)$, but
in  the last one  it denotes the action of  $(F(K),\succ, U(\mathfrak{g}))$ and
of  $( F(L), \succ, U(\mathfrak{l}))$, respectively.

{\bf (3)} When the   regular module
 $(H^*,\prec, H)$ is taken in consideration, the  following chains of 
equalities determine  the  action of  $k'$ and $\lambda '$, respectively:
\begin{equation}\begin{array}{lll}
  \langle (\kappa l)\prec k', k \lambda \rangle &= & 
       \langle \kappa l, k'\succ (k \lambda) \rangle=
       \langle \kappa l, k' k \lambda \rangle \\[0.2cm]
  &= &   \langle \kappa , k' k \rangle \langle l, \lambda \rangle \quad =
  \langle \kappa \prec k', k \rangle \langle l, \lambda \rangle 
= \langle (\kappa \prec k')l, k \lambda \rangle \ ;
   \end{array} \end{equation}
\begin{equation} \begin{array}{lll}\label{lambdaaction}
  \langle (\kappa l)\prec \lambda', k \lambda \rangle &= & 
  \langle \kappa l, \lambda'\succ( k \lambda) \rangle =    
  \langle \kappa l,  k (\lambda' \ract k) \lambda \rangle     \\[0.2cm]
 &= &  \langle \kappa, k  \rangle \langle l,  (\lambda' \ract k) \lambda
 \rangle= \langle \kappa, k  \rangle \langle l,  (\lambda' \ract k)\rangle
       \langle l,  \lambda \rangle \\[0.2cm]
 &= &  \langle \kappa, k  \rangle
  \langle l^{(1)}, \lambda' \rangle \langle l^{(2)},  k \rangle
       \langle l,  \lambda \rangle=  \langle \kappa l^{(2)}, k  \rangle
  \langle l^{(1)}, \lambda' \rangle
   \langle l,  \lambda \rangle  \\[0.2cm]
~& = &  \langle \kappa l^{(2)}l , k \lambda  \rangle
  \langle l^{(1)}, \lambda' \rangle =
   \langle \kappa \langle l^{(1)}, \lambda' \rangle l^{(2)}l , k \lambda  \rangle \ .
    \end{array}\end{equation}
\end{proof}
Note that: (i) the action (\ref{lambdaaction}) is described in terms of  the  structure
of  $U(\mathfrak{l})$ as right $F(K)$--comodule; and (ii)  except the term
$\langle l^{(1)}, \lambda' \rangle l^{(2)}$ including a coaction,  the  action  on the
regular modules appears described by means of  other actions,  most of  them regular.

From a computational point of  view  the  following
proposition and its   corollary are very useful, since they allow to reduce  the 
description of  the regular modules to the study of  the  action of  $K$
on $L$ associated to  the  structure of 
$U(\mathfrak{k})$--module of  $F(L)$. 

Let us start fixing   the 
notation to be used. Let $r$ and $s$ be the dimensions
of  the groups $K$ and $L$, respectively. Let us consider the basis $(k_i)_{i=1}^r$ of 
$\mathfrak{k}$ and $(l_j)_{j=1}^s$ of $\mathfrak{l}$, and the local coordinate systems of
second kind associated to the above bases $(\kappa_i)_{i=1}^r$ and
$(\lambda_j)_{j=1}^s$.
Remember that using  multi-index notation one has
 \begin{equation} \label{bases}
    \langle k_n, \kappa^{n'} \rangle = n! \; \delta^{n'}_n,  \qquad
    \langle l_m, \lambda^{m'} \rangle = m! \; \delta^{m'}_{m},\qquad
n, n' \in \mathbb{N}^r,\ m, m' \in \mathbb{N}^s.
\end{equation}
Finally, let us denote by $k$ the  inverse map of the coordinate
system $(\kappa_i)$, i.e., 
\begin{equation}
    \begin{array}{cccc} k: & \mathbb{R}^r & \longrightarrow  & K \\
                   &  t& \mapsto & e^{t_1 k_1} e^{t_2 k_2} \cdots e^{t_r k_r}.
    \end{array} \end{equation}

\begin{proposition} 
For every  $\lambda \in F(L)$ and  $l\in L$ the following relation holds
  \begin{equation}\langle l^{(1)}, \lambda \rangle l^{(2)} =
      \sum_{n \in \mathbb{N}^r} \frac{1}{n!} 
     \langle k_n \lact l, \lambda \rangle \kappa^n
         = \lambda(k(\kappa)\lact l).  
 \end{equation}
\end{proposition}
\begin{proof}[Proof]
Let us rewrite the  coaction  at  the  right of  $F(K)$ on $l \in L$ as
 \begin{equation} \label{crcoaccion}
   \rcact l\equiv  l^{(1)}\otimes l^{(2)} =\sum_{(m,n)\in \mathbb{N}^s \times
\mathbb{N}^r}
           [l]^m_n \, l_m \otimes \kappa^n.
 \end{equation} 
The   pairing defined in the   bicrossproduct in accordance with 
Theorem \ref{tbasesduales} allows to obtain the coordinates of  $\rcact l$
in terms of  the  action  (dual of  the  coaction) of  $U(\mathfrak{k})$
on $U(\mathfrak{l})$
 \begin{equation}
    [l]_m^n=
       \frac{1}{  m! n!} \langle \lambda^m \otimes k_n, \rcact l
\rangle =\frac{1}{  m! n!} \langle \lambda^m, k_n \lact l \rangle. 
 \end{equation}
Inserting the last expression in (\ref{crcoaccion}) 
one easily gets
\begin{equation} \label{otraexpr}
  \begin{array}{lll}
    \langle l^{(1)}, \lambda \rangle l^{(2)} &= &
      \sum_{(m,n) \in \mathbb{N}^s \times \mathbb{N}^r}
       \frac{1}{m!n!} 
     \langle \lambda^m, k_n \lact l \rangle
     \langle l_m, \lambda \rangle \kappa^n  \\[0.3cm]
   &= &    
      \sum_{(m,n) \in \mathbb{N}^s \times \mathbb{N}^r}
       \frac{1}{m!n!} 
     \langle k_n \lact l,  \lambda^m \rangle
     \langle l_m, \lambda \rangle \kappa^n.
    \end{array}
 \end{equation}
The  sum on $m$ gives account of  the  action of  the  $T$--matrix
associated to the  pair $(U(\mathfrak{l}), F(L))$ and, hence, 
 the  expression (\ref{otraexpr}) is
simplified getting
\begin{equation}\langle l^{(1)}, \lambda \rangle l^{(2)} =
      \sum_{n \in \mathbb{N}^r} \frac{1}{n!} 
     \langle k_n \lact l, \lambda \rangle \kappa^n.\end{equation}
On the other hand, since
\begin{equation}
    \lambda(k(t)\lact l)= \langle k(t) \lact l, \lambda \rangle=
    \sum_{n\in \mathbb{N}^r}
     \frac{1}{n!} \langle k_n \lact l, \lambda \rangle t^n,
 \end{equation}
in order to get  $\langle l^{(1)}, \lambda \rangle l^{(2)}$
it is enough to perform the   formal  substitution
$t_i \to \kappa_i$
in  the  expression $\lambda(k(t)\lact l)$.
\end{proof}

 The  above proposition allows to give an expression for
$\langle l^{(1)}, \lambda \rangle l^{(2)}$
completely independent of  the bases chosen in the algebras.

\begin{corollary} Let $ \hat{l}$ be the  map 
\begin{equation} \label{proyeccion}
   \begin{array}{cccc}
     \hat{l}: & K & \longrightarrow & L \\
              & k & \mapsto & k \lact l
   \end{array}
\end{equation}
 projecting the  group $K$ on the  orbit passing through $l \in L$. Then, for any
$\lambda \in F(L)$ and any $l\in L$ one has
  \begin{equation}
\langle l^{(1)}, \lambda \rangle l^{(2)} = \lambda \circ \hat{l}. 
\end{equation}
\end{corollary}

Taking into account that in
$(U(\mathfrak{l}),\succ, F(L))$ it is verified that
\begin{equation}
 \lambda \succ l= \lambda(l) l,
  \qquad \forall \lambda \in F(L), \; \forall l \in L
\end{equation}
Theorem \ref{crbicross0} can be rewritten in a more explicit way.

\begin{theorem} \label{crbicross}
The  action  on each of  the four regular $H$--modules is:
\begin{equation}\begin{array}{lllll} 
(H,\prec, H)&:&\ \ (k\lambda) \prec k'=  kk'(\lambda \ract k') , 
\qquad\quad  & (k\lambda) \prec \lambda'=  k\lambda \lambda' ;\\[0.2cm]
 (H^*,\succ, H)&:&\ \  k' \succ (\kappa l)  =  (k'\succ \kappa) (k'\lact l),
\qquad\quad  &\lambda' \succ (\kappa l)  = \lambda'( l)\kappa  l ;\\[0.2cm]
(H,\succ, H)&:&\ \  k' \succ (k\lambda) =  k'k\lambda ,
\qquad\quad  &\lambda' \succ (k\lambda) =  k (\lambda' \ract k) \lambda ;\\[0.2cm]
(H^*,\prec, H)&:&\ \   (\kappa l) \prec k'  =  (\kappa \prec k') l ,
\qquad\quad  &(\kappa l) \prec \lambda'  =\kappa (\lambda'\circ \hat{l}) l ;
\end{array}
\end{equation}
where
  $k,k' \in K$, $\lambda, \lambda' \in F(L)$, $\kappa \in F(K)$ and $l \in L$.
\end{theorem} 
The   result of this  theorem does not make reference to  the 
nature of   Lie groups  $K$ and $L$, since it is formulated in terms of the
regular actions and associated
ones to the bicrossproduct structure. Thus, the theorem may be applied to other
kinds of groups.

Note that, in general,  the  action of  $K$ on $L$ is not
globally defined. Hence, $ \hat{l}$ (\ref{proyeccion})
 only projects, in reality,  a neighbourhood of  the  identity into the orbit of $l$.
Henceforth, $\lambda \circ \hat{l}$
does not define, in general, a map over the whole $K$ and  the  expression
$(\kappa l) \prec \lambda'  = \kappa (\lambda'\circ \hat{l}) l $ only has
sense enlarging  the   space
$F(K)$, for instance, including it inside  spaces of  formal series.

As a  conclusion, we can say that in the description of the regular actions the  
computation  of  the left  action of the  group $K$ on the   group $L$ is really
the most important fact. From this point of view, the deformations used in this
work may be interpreted as one-parameter families of nonlinear actions 
homotopically equivalent to the linear actions of the nondeformed cases.
\subsect{Co-spaces and induction}

In the context of noncommutative geometry the  
manifold $X$ is  replaced  by the algebra ${F}(X)$ of
${\cal C}^\infty$ $\C$--valued functions on $X$  as well as the Lie group $G$ by the
enveloping  algebra $U(\mathfrak{g})$ of its Lie algebra $\mathfrak{g}$. Since  $({F}(X),
\triangleright, U(\mathfrak{g}))$ is a module algebra over the Hopf algebra
$U(\mathfrak{g})$, we can generalize the concept of $G$--space in algebraic terms \cite{ao00}.

Let $H$ be a Hopf algebra. A left (right) $H$--co-space is a module algebra
$(A, \triangleright, H)$ ($(A, \triangleleft, H)$).

The morphisms among \mbox{$H$--co-spaces} are the morphisms of
$H$--modules and the concepts of subco-space or quotient co-space are
equivalent to  module subalge\-bra  or quotient module algebra, respectively.
We have adopted the term of {\em co-space} instead of {\em space} to stress the dual character of $A$ as
way of describing the initial geometric object.

Given a pair of algebras with a
non-degenerate pairing $(H,  H', \langle \cdot, \cdot \rangle)$, we obtain, via 
dualization of the regular actions, the regular $H$--co-spaces
$(H',\succ, H)$ and $(H',\prec,H)$. 
 
The  explicit description  of  the four regular modules studied in the previous subsection  allows a
complete  analysis of the  representations of the algebra $H= U(\mathfrak{k}) \RL F(L)$
induced by the one-dimensional modules of  the  commutative sector. As we will see, the left co-space
$(H^*,\succ, H)$ characterizes the carrier space of the induced representation and the right co-space
$(H^*,\prec,H)$ determines the induced action of $H$ on the carrier space. 

Firstly,  remember that the set of  characters of the algebra $F(L)$ is its spectrum. An important
theorem by Gelfand and Naimark \cite{Gel43} establishes the following isomorphism
  \begin{equation}
        \text{Spectrum} \; F(L) \simeq L.
 \end{equation}
 Fixed $l\in L$, the  character    (or the  corresponding right $F(L)$--module
over $\mathbb{C}$) is given by
\begin{equation}  \label{car}
   1 \dashv \lambda = \lambda(l), \quad \lambda \in F(L).
\end{equation}
In order to construct  the  representation of  $H= U(\mathfrak{k}) \RL F(L)$
induced by (\ref{car}) let us start determining the carrier space 
$\mathbb{C}^\uparrow \subset H^*$. The   element $f \in H^*$ satisfies  the 
equivariance condition if it verifies 
 \begin{equation}
    \lambda \succ f = \lambda(l) f, \qquad \forall \lambda \in F(L).
 \end{equation}
Expanding $f$ in terms of the  bases  of $\mathfrak{k}$ and $\mathfrak{l}$
\begin{equation}
 f= \sum_{(m,n)\in \mathbb{N}^r \times \mathbb{N}^s} f_m^n \kappa^m l_n ,
\end{equation}
the equivariance condition gives the following relation among the
coefficients $f^n_m$
\begin{equation}
f^n_m= \frac{1}{m! n!}  f^0_m \lambda^n(l),
\qquad m \in \mathbb{N}^r, n \in \mathbb{N}^s.
\end{equation}
Hence, the general  solution is:
\begin{equation}\label{solucioncovariancia}
f = \left( \sum_{m \in \mathbb{N}^r} \frac{1}{m!} f_m^0 \kappa^m \right)
    \left( \sum_{n \in \mathbb{N}^s } \frac{1}{n!} \lambda^n(l) l_n \right),
  \qquad f_m^0 \in \mathbb{C}.
\end{equation}
Taking into account the  definition of the
second kind coordinates $\lambda_j$ over the   group $L$,  the  expression
(\ref{solucioncovariancia}) can be rewritten in a more compact form
\begin{equation}
     f= \kappa l, \qquad  \kappa \in F(K).
\end{equation}
In other words, the   carrier space  of  the  induced
representation admits a natural description
in terms of  products  {\em function/element}, introduced
in (\ref{tablita}), instead of terms of monomial bases.

The  right regular action  describes  
the  action  on the  induced module, which  can be translated to
$F(K)$ using the isomorphism
$ F(K)  \rightarrow  \mathbb{C}^\uparrow \ (  \kappa  \mapsto  \kappa l)$:
\begin{equation}
 \kappa \dashv k= \kappa \prec k, \qquad
 \kappa  \dashv \lambda=   \kappa  (\lambda \circ \hat{l}).
\end{equation}
Comparing these expressions with those of Theorem \ref{tind1} we observe 
that  the  action of  the  subalgebra $U(\mathfrak{k})$ is given by 
the  regular action.  The  action of  the  subalgebra $F(L)$ is of 
 multiplicative kind and  the  evaluation of the corresponding  factor, from a
computational point of view, is essentially reduced to  obtain  the one-parameter
flows associated to  the  action of  $K$ on $L$ derived of  the bicrossproduct 
structure  of the algebra $H$.

\subsect {Equivalence and unitarity of  the induced representations}

Let $\dashv_l$ be the  representation of  $H= U(\mathfrak{k})
\RL F(L)$ induced by $l\in L$ and  $f_k$ the automorphism of  $F(K)$ 
given by  the  regular action of  an element $k\in K$, i.e.,
$     f_k(\kappa)= k \succ \kappa$. Since
$$
[k \succ (\lambda \circ \hat{l})](k')=
( \lambda \circ \hat{l})(k'k)= \lambda((k'k) \lact l)=
[\lambda \circ \widehat{k \lact l}](k') ,
$$
then one has that
$$  
f_k(\kappa \dashv_l \lambda)=  k \succ [\kappa  (\lambda \circ \hat{l})]=
  (k \succ \kappa )[k \succ (\lambda \circ \hat{l})]=
     f_k(\kappa) [ \lambda \circ \widehat{k\lact l}]=
   f_k(\kappa) \dashv_{\widehat{k\lact l}} \lambda. 
$$
Taking into account, besides, that  the  action
of the  subalgebra $U(\mathfrak{k})$ on the  induced module   is not affected  by 
the  choice of the element $l$ in $L$, we conclude that the $H$--modules
$(\mathbb{C}^\uparrow, \dashv_l, H )$ and 
$(\mathbb{C}^\uparrow, \dashv_{\widehat{k\lact l}} , H )$ are isomorphic via $f_k$. 

 The   problem of  the  unitarity of 
the induced representation  passes, firstly, for choosing a $*$--structure in
$H$.  The  usual determination is to consider `hermitian operators' a  family of 
generators of  $H$, but troubles, related with the real or complex nature of the
deformation parameter, may appear \cite{bgst98,Bon98}. The   point of  view adopted here
allows  a simple solution of the problem:
$U(\mathfrak{k})$ and $ F(L)$ 
carry associated $*$--structures in a natural way. Explicitly, 
\begin{equation} \label{involus}
   \begin{split}
       k^*= & \ k^{-1},  \qquad \forall k \in K,  \\[0.2cm]
     \lambda^*(l)= &\  \overline{\lambda(l)},
                 \qquad \forall \lambda \in F(L), \ \forall l \in L.
\end{split}
 \end{equation}
Choosing in $H$  the  $*$--structure associated to those given by
(\ref{involus}), according to Theorem \ref{estructuraestrellabicross},
the  problem of  the  unitarization is easily solved.
Firstly,  the  action of  the elements $k \in K$ shows that the  
space $F(K)$ has to be restricted to the  
square-integrable functions  with respect to   the right invariant Haar
measure $\mu$ over $K$ (i.e., 
$\mu( k \lact A) = \mu(A)$
with $ A$ a $\mu$--measurable set in $K$). In fact, it is necessary to restrict
the  space ${\cal H}=L^2(K,\mu)$ and to consider only the space ${\cal H}_\infty$ of 
${\cal C}^\infty$ functions, since  the  Lie algebra
$U(\mathfrak{k})$ acts  by means  of differential  operators over these functions. On
the other hand, the elements of  $F(L)$ act by  a multiplicative factor and
impose a new restriction in ${\cal H}_\infty$ because only the  functions
$\kappa$ such that  $\kappa (\lambda\circ \hat{l})$
is also square-integrable  (supposing that  the  action  is
global in  the orbit of $l$) will be admissible.
If $K$ is compact all that is automatically verified and in the opposite case
there is a condition over the   vanishing order  of 
$\kappa$ at the infinity points.
The results of  this  discussion are summarized in the   following theorem.
\begin{theorem}  \label{repequiv}
Let us consider an element  $l\in L$ supporting a global  action of the  group
$K$.  The carrier space, $\mathbb{C}^\uparrow$, of the  representation of  $H$
induced by the  character    determined by $l$ is the  set of  elements of  $H^*$
of  the  form
 \begin{equation} \kappa l, \qquad \kappa \in F(K).
\end{equation}
There is an isomorphism between
$\mathbb{C}^\uparrow$ and $F(K)$ given by the  map $\kappa \mapsto \kappa l$.
The  action  induced  by the elements of  the  form 
$k \in K$ and $\lambda \in F(L)$ in the  space $F(K)$ is
\begin{equation}
\begin{array}{lll}
        \kappa \dashv k &= & \kappa \prec k \\[0.2cm]
        \kappa  \dashv \lambda &= &
             \kappa  (\lambda \circ \hat{l}) .
\end{array}
\end{equation} 
The modules induced by $l$ and $k\lact l$ are isomorphic. So,
the induction algorithm establishes a correspondence between the  
space of  orbits $L/K$ and the set of  equivalence classes of  representations.

If the   group $K$ is compact  the  induced representation is unitary in
the  space $L^2(K)$, of  square-integrable functions with respect to  the right
invariant Haar measure, when the  $*$--structure given by  Theorem 
\ref{estructuraestrellabicross}, applied to the natural structures  of 
the factors of the bicrossproduct $H= U(\mathfrak{k})\RL F(L)$, is considered.
\end{theorem}
\subsect{Local  representations}
\medskip

The called local representations \cite{olmo84} in the deformed version appear when one 
induces from representations of  the  subalgebra $U(\mathfrak{k})$. Let us consider
the following  character of  $U(\mathfrak{k})$
  \begin{equation} 
\kappa \in \text{Spectrum} \, U(\mathfrak{k}) \subset F(K), \qquad
      k \vdash 1 = \kappa(k).
\end{equation}
Since the  Hopf algebra $U(\mathfrak{k})$ is, in general, non commmutative, the set of
characters may be very reduced, even it may be   generated  only by the  counit. For
this reason an interesting problem to be researched in the future is the  
study of the representations  induced by representations of  $U(\mathfrak{k})$ of 
dimension greater that one.

The  carrier space of  the  representation induced by $\kappa$ is determined by
the following equivariance condition
 \begin{equation}
    f \prec k=  \kappa(k) f, \qquad  \forall k \in U(\mathfrak{k}).
 \end{equation} 
The algebra $H^*$ can be consider as a  left free $F(K)$--module and, hence, it is
possible to fix a basis $(l_j)_{j\in J}$ of 
$U(\mathfrak{l})$ such that $f \in H^*$ can be expressed in a unique form as
 \begin{equation} \label{freemodule}
     f= \sum_{j\in J} \kappa_j l_j, \qquad \kappa_j \in F(K).
 \end{equation}
The equivariance condition can be written now as 
\begin{equation}
  \sum_{j \in J} (\kappa_j \prec k) l_j = \sum_{j \in J} \kappa(k) \kappa_j l_j,
\qquad
        \forall k \in U(\mathfrak{k}).
\end{equation}
Taking into account that the elements $l_j$ constitute a basis of the $F(K)$--module $H^*$ the corresponding
coefficients can be equating, obtaining
\begin{equation}\label{previousequality}
  \kappa_j \prec k = \kappa(k) \kappa_j, \qquad  \forall k \in U(\mathfrak{k}). 
\end{equation}
The  previous equality (\ref{previousequality}) implies that
\begin{equation}
\kappa_j(kk')= \kappa(k) \kappa_j(k'), \qquad \forall k, k' \in K.
\end{equation}
Choosing $k'$ equal to the identity  element $e \in K$, one gets
\begin{equation}
  \kappa_j = \kappa_j(e) \kappa, \qquad \forall j \in J,
\end{equation}
in this way the elements of  $\mathbb{C}^\uparrow$ are of the form
\begin{equation}
         f = \kappa l, \qquad  l \in U(\mathfrak{l}).
\end{equation}
 The  map $U(\mathfrak{l})\to \mathbb{C}^\uparrow$, defined by $l \mapsto \kappa l$,
is an  isomorphism of vector spaces.
The  representation   can be  realized in this way  over $U(\mathfrak{l})$ and
the  final  result is
\begin{equation} \label{reploc}
 \begin{array}{c}
        k \vdash l= \kappa(k) \; k \lact l ,\\[2mm]
        \lambda \vdash l = \lambda(l) l .
   \end{array}
\end{equation}

\sect{Examples}\label{ejemplos}

\subsect{Null-plane quantum Poincar\'e algebra}

The null-plane quantum deformation of the $(1+1)$ Poincar\'e algebra,
$U_{z}(\mathfrak{p}(1,1))$, is a $q$--deformed Hopf algebra that in a null-plane basis,
$\{ P_+,P_-, K\}$, has the form \cite{Bal95f,Bal95g} 
\begin{equation}
  \begin{array}{c}
[K, P_+]= \frac{-1}{z}(e^{-2 z P_+} - 1), \qquad
[K, P_-]= -2 P_-, \qquad [P_+, P_-]= 0; \\[0.2cm]
\Delta P_+= P_+ \otimes 1 + 1 \otimes P_+, \qquad
\Delta X = X \otimes 1 + e^{-2 z P_+} \otimes X , \quad X \in \{P_-, K\};\\[0.2cm]
 \epsilon(X)= 0, \quad  X \in \{P_\pm, K\}; \\[0.2cm]
 S(P_+)=- P_+, \qquad S(X)=- e^{2 z P_+} X, \qquad X \in \{P_-, K\}.
  \end{array}
\end{equation}
It has also the structure of bicrossproduct \cite{PerJ96}
$$
U_z (\mathfrak{p}(1,1))= {\cal K} \RL {\cal L},
$$ 
where $\cal K$ is a commutative and
cocommutative Hopf algebra generated by $K$, and $\cal L$ is the commutative Hopf subalgebra
of $U_z (\mathfrak{p}(1,1))$ generated by $P_+$ and $P_-$.

The right action of $\cal K$ on $\cal L$ comes determined by
\begin{equation}
 P_+ \ract K= \frac{1}{z}(e^{-2 z P_+}-1), \qquad
 P_- \ract K= 2 P_- .
\end{equation}
The left  coaction of $\cal L$ over the generator of $\cal K$ is
\begin{equation}
  K \lcact = e^{-2 z P_+} \otimes K.
\end{equation}

In the dual Hopf algebra  $F_z (P(1,1))= {\cal K}^* \LR {\cal L}^*$
let us  denote by $\varphi$ the generator of ${\cal K}^*$ and
by $a_\pm$ those of ${\cal L}^*$. The left action of ${\cal L}^*$ on ${\cal K}^*$ is given
by
 \begin{equation} 
a_+ \lact \varphi=  2 z (e^{-\varphi} -1),
    \qquad a_- \lact \varphi= 0, 
\end{equation}
and the right coaction of ${\cal K}^*$ over the generators of ${\cal L}^*$ by
 \begin{equation} 
\rcact a_\pm= a_\pm \otimes e^{\mp 2 \varphi}.  
\end{equation}
With these actions we obtain the Hopf algebra structure of $F_z (P(1,1))$
\begin{equation}
\begin{array}{c}
 [a_+, a_-]= -2 z a_-, \qquad [a_+, \varphi]= 2 z (e^{-\varphi} -1), 
\qquad [a_-, \varphi]=0;\\[0.2cm]
     \Delta a_\pm = a_\pm \otimes e^{\mp 2 \varphi} + 1 \otimes a_\pm \quad, \qquad
  \Delta \varphi = \varphi \otimes 1 + 1 \otimes \varphi; \\[0.2cm] 
\epsilon(f)= 0, \quad  f \in \{a_\pm, \varphi \}; \\[0.2cm] 
S(a_\pm)= - a_\pm e^{\pm \varphi}, \qquad S(\varphi )=- \varphi.
\end{array}
\end{equation}

Theorem \ref{tbasesduales} allows to obtain easily a pair of dual bases in such a
way that the duality between $U_z (\mathfrak{p}(1,1))$ and $F_z (P(1,1))$
is explicitly given by the  pairing
\begin{equation}
   \langle K^m P_-^n P_+^p, \varphi^q a_-^r a_+^s \rangle =
   m! n! p! \; \delta^m_q \delta^n_r \delta^p_s.
\end{equation}

Now let us consider the  bicrossproduct structure of $U_{z}(\mathfrak{p}(1,1))$ as
follows
 \begin{equation}
    U_{z}(\mathfrak{p}(1,1))=
     U(\mathfrak{k})  \RL F(T_{{z}, 2}),
\end{equation}
where $\mathfrak{k}$ is the one-dimensional Lie algebra generated by $K$
 and  the group  $T_{{z}, 2}$ is  a  deformation of the additive group
$\mathbb{R}^2$ defined by the law
\begin{equation} 
(\alpha'_-,\alpha'_+ )(\alpha_-,\alpha_+ )=
    (\alpha'_- e^{-2{z} \alpha'_+} \alpha_-, \alpha'_+ + \alpha_+).
\end{equation}
The functions
\begin{equation} 
P_-(\alpha_-,\alpha_+ )= \alpha_-, \quad
   P_+(\alpha_-,\alpha_+ )= \alpha_+
\end{equation}
define a global  chart on $T_{{z}, 2}$. The   
$U(\mathfrak{k})$--module algebra structure of  $F(T_{{z},2})$ taking part in the
bicrossproduct is given by 
\begin{equation}
   P_- \ract K= 2 P_-, \qquad
   P_+ \ract K= \frac{1}{{z}}(e^{-2 {z} P_+}-1).
\end{equation} 
Hence, the   vector field associated to $K$ is
\begin{equation} 
\hat{K}= 2 P_- \frac{\partial}{\partial P_-}+
           \frac{1}{{z}}(e^{-2 {z} P_+}-1) \frac{\partial}{\partial P_+}.
\end{equation}

\subsubsect{One-parameter flow}

The   vector field $\hat{K}$ has a unique  equilibrium point at $(0,0)$, which
has hyperbolic nature. The function 
\begin{equation}
    h= P_-(e^{-2{z} P_+}-1)
\end{equation}
is a first integral  of  $\hat{K}$. The   computation of  the integral curves require to 
solve the   differential system 
\begin{equation}
     \dot{\alpha}_-=  2 \alpha_-, \qquad
     \dot{\alpha}_+=  \frac{1}{z} (e^{-2{z} \alpha_+}-1).
\end{equation}
If $z>0$ the integral curves placed in  the  region $\alpha_+ <0$ are given by 
\begin{equation} 
    \alpha_-(s)=  c_1 e^{2s}, \qquad
    \alpha_+(s)= \frac{1}{2{z}} \ln(1-e^{-2(s-c_2)}).
 \end{equation} 
The second order system associated to them is 
 \begin{equation}
      \ddot{\alpha}_-(s)=  4 {\alpha}_-(s), \qquad
      \ddot{\alpha}_+(s)=  -\frac{2}{{z}} e^{-2{z} {\alpha}_+}(e^{-2 {z}{\alpha}_+}-1).
       \end{equation}
These equations may be interpreted as particles moving  over  a straight line under  the 
action of repulsive potentials.  From  the  expression of  the integral curves we get the
following flow
\begin{equation}
    \Phi^s(\alpha_-, \alpha_+)= (\alpha_- e^{2s},
   \frac{1}{2{z}} \ln(1-e^{-2s}(1-e^{2{z} \alpha_+}))).
      \end{equation}
If, for example, we suppose that ${z} > 0$ then  the  curve that starts at the point
$(\alpha_-, \alpha_+)$ is defined in the  interval
\begin{equation} s \in \left\{ \begin{array}{lcc}
       (\frac{1}{2} \ln(1-e^{2{z} \alpha_+}), + \infty) & & \alpha_+ < 0 ,\\[2mm]
       (-\infty, +\infty) & & \alpha_+ \geq 0 . \end{array} \right.
\end{equation}
Hence,  the  expression
\begin{equation}
     e^{sK} \lact (\alpha_-, \alpha_+)= (\alpha_- e^{2s},
   \frac{1}{2{z}} \ln(1-e^{-2s}(1-e^{2{z} \alpha_+})))
\end{equation}
defines a local  action (except in the nondeformed  limit ${z} \rightarrow 0$, where
the  action  is global) of  $\mathfrak{K}$ (the Lie group associated to the Lie algebra $\mathfrak{k}$)
on $T_{{z}, 2}$. 
The  action  decomposes  $T_{{z}, 2}$ in three strata:

\ \ i) the   point at the origin, whose isotropy  group is $\mathfrak{K}$,

\ ii) the four orbits constituted by the  semiaxes,

iii) the   rest of the set $T_{{z}, 2}$. This last stratum has a foliation by
one-dimensional orbits: deformed hyperbolic branches.

\subsubsect{Regular co-spaces}\label{co-spaces1}
The elements of $F_{z}(P(1,1)$ can be written as
\begin{equation}
\phi(\alpha_-, \alpha_+), \qquad  \phi \in F_{z}(\mathfrak{K}), \ \ (\alpha_-, \alpha_+) 
\in T_{{z}, 2} 
 \end{equation}
instead of the monomials $\varphi^q a_-^r a_+^s$. The  expression $\phi(\alpha_-, \alpha_+)$ does not
denote a  function, $\phi$, at the point
$(\alpha_-, \alpha_+)$ but the   product of these two elements in the   algebra
$F_{z}(P(1,1))$.

  The  structure of the  regular co-space
 $(F_{z}(P(1,1)), \prec, U_{z}(\mathfrak{p}(1,1)))$ is immediately obtained using 
   Theorem~\ref{crbicross}.  So,
\begin{equation}\begin{split}
 (\phi(\alpha_-, \alpha_+))\prec e^{s K}= & \phi(e^{s K} \, \cdot \, )
       (\alpha_-, \alpha_+), \\[2mm] 
(\phi(\alpha_-, \alpha_+))\prec P_-= & \phi \alpha_- e^{2\varphi}
       (\alpha_-, \alpha_+), \\[2mm]
 (\phi(\alpha_-, \alpha_+))\prec P_+= & \phi \frac{1}{2{z}}
   \ln(1-e^{-2\varphi}(1-  e^{2 {z} \alpha_+}))  (\alpha_-, \alpha_+),
  \end{split}
 \end{equation}
The  action  on
 $(F_{z}(P(1,1)), \succ, U_{z}(\mathfrak{p}(1,1)))$ is given by 
\begin{equation}\label{rpz2}\begin{split}
  e^{s K} \succ (\phi(\alpha_-, \alpha_+))= & \phi(\, \cdot \, e^{s K})
       (\alpha_- e^{2s},\frac{1}{2{z}}
   \ln(1-e^{-2s}(1-  e^{2 {z} \alpha_+}))), \\[2mm]
 P_- \succ (\phi(\alpha_-, \alpha_+))= & \alpha_-  \phi
       (\alpha_-, \alpha_+), \\[2mm]
 P_+ \succ (\phi(\alpha_-, \alpha_+))= & \alpha_+  \phi
       (\alpha_-, \alpha_+).       
  \end{split}
 \end{equation}
In the above expressions the dot stands for the argument of the function $\phi=\phi(\; \cdot\; )$, and
$\varphi$ denotes  the  natural coordinate function over  the group $\mathfrak{K}$.

Note that the elements $(\alpha_-, \alpha_+) \in T_{{z},2}$ describe  the  subalgebra of 
$F_{z}(P(1,1))$ generated by $a_-$ and $a_+$. The pair $(\alpha_-, \alpha_+)$ is an eigenvector of 
the endomorphisms associated to  the  action  (\ref{rpz2}) of  the generators
$P_-$ and $P_+$. This fact, together with 
 the  action of  $\mathfrak{K}$ on $T_{{z},2}$,
guarantees  that the  subalgebra generated by $a_-$ and $a_+$ is stable under  the 
action  (\ref{rpz2}).

\subsubsect{Induced representations}
 The  representation of  $U_z(\mathfrak{p}(1,1))$
 induced  by the  character    $(\alpha_-, \alpha_+) \in T_{{z},2}$ is given
according to   Theorem~\ref{repequiv} by the following expressions
\begin{equation}\begin{split}
 \phi \dashv  K= & \phi', \\[2mm]
 \phi \dashv  P_-= & \phi \alpha_- e^{2\varphi}, \\[2mm]
 \phi \dashv P_+= & \phi \frac{1}{2{z}}
   \ln(1-e^{-2\varphi}(1-  e^{2 {z} \alpha_+})).
  \end{split}
 \end{equation}
Choosing a representative in  each orbit one gets a 
representative  of  every equivalence classes of induced  representations.  For
instance, the  representation induced by the equilibrium point
$(0, 0) \in T_{{z},2}$ is
\begin{equation}\begin{split}
 \phi \dashv  K= & \phi', \\[2mm]
 \phi \dashv  P_\mp= & 0.
  \end{split}
 \end{equation}

The local  representations induced by the  character   
 \begin{equation} 
K^m \vdash 1= c^m,
\end{equation} 
of  the  subalgebra $U(\mathfrak{k})$ are given, according to  
(\ref{reploc}), by
\begin{equation} \begin{split}
     e^{s K} \vdash (\alpha_-, \alpha_+)= &
         e^{s c}   (\alpha_- e^{2 s},
   \frac{1}{{z}} \ln(1-e^{-2 s}(1-e^{2{z} \alpha_+}))), \\[2mm]
     P_- \vdash (\alpha_-, \alpha_+)= & \alpha_- (\alpha_-, \alpha_+), 
\\[2mm]
     P_+ \vdash (\alpha_-, \alpha_+)= & \alpha_+ (\alpha_-, \alpha_+).
   \end{split}\end{equation}

\subsect{Non-standard quantum  Galilei algebra}
\label{deformacionnoestandarg}

The non-standard quantum Galilei algebra $U_z(\mathfrak{g}(1,1))$ is isomorphic to the
quantum Heisenberg algebra $H_q(1)$ \cite{cgst92,cgst91} and to the deformed
Heisenberg--Weyl algebra
$U_\rho(HW)$ \cite{apb96}. It can be obtained by contraction \cite{apb96} of a
non-standard deformation of the Poincar\'e algebra \cite{Bal95g} (the null-plane
quantum Poincar\'e).

The deformed Hopf algebra $U_z(\mathfrak{g}(1,1))$ has the following  structure
\begin{equation}\label{estrgalz}
\begin{array}{c}
[H, K]= - \frac{1- e^{-4 z P}}{4z}, \qquad  [P, K]= 0, \qquad [H, P]=0;  \\[2mm]
\Delta P = P \otimes 1 + 1 \otimes P,  \qquad
\Delta X = X \otimes 1 + e^{- 2 z P} \otimes X, \quad X \in \{H, K \}; \\[2mm]
 \epsilon(X)= 0,  \qquad  X \in \{H, P, K\}; \\[2mm]
S(P)= -P,  \qquad  S(X)= - e^{2 z P} X,   \quad X \in \{H, K\}.
  \end{array}
\end{equation}

In  \cite{PerJ96} it was proved that  $U_{z}(\mathfrak{g}(1,1))$ has
structure of bicrossproduct 
$$
U_z(\mathfrak{g}(1,1))={\cal K}\RL {\cal L},
$$ 
where ${\cal L}$ is the commutative and non-cocommutative Hopf subalgebra 
$U_z (\mathfrak{t}_2)$ generated by $P$ and $H$, and $\cal K$ is the commutative
and cocommutative Hopf algebra (it is not a Hopf subalgebra of
 $U_z (\mathfrak{g}(1,1))$) generated by  $K$.

The right  action of ${\cal K}$ on ${\cal L}$ is given by
 \begin{equation}
     P \ract K= [P,K]= 0, \qquad H \ract K= [H,K]=
 - \frac{1- e^{-4z P}}{4 z}.
 \end{equation}
The left coaction of ${\cal L}$ over the generator of ${\cal K}$ is
 \begin{equation} 
K \lcact= e^{-2 z P} \otimes K . 
\end{equation}

The corresponding function algebra $F_z (G(1,1))$
has a  bicrossproduct structure dual of the above one 
$$
F_z (G(1,1))= {\cal K}^* \LR {\cal L}^* .
$$
Let $v,x$ and $t$ be the  generators
dual of $K,P$ and $H$. The  action of ${\cal L}^*$ on 
${\cal K}^*$ is
 \begin{equation}
         x \lact v= - 2 z v, \qquad t \lact v= 0,
 \end{equation}
and the  coaction of ${\cal K}^*$ over the generators of  ${\cal L}^*$ is
 \begin{equation}
         x \lcact = 1 \otimes x, \qquad t \lcact = 1 \otimes t.
 \end{equation}
Action and coaction allow to obtain the  Hopf algebra structure of $F_z (G(1,1))$
\begin{equation}
\begin{array}{c}
[t, v]= 0,  \qquad [x, v]= - 2 z v,  \qquad [t,x] =2 z t \quad;  \\[2mm]
\Delta t = t \otimes 1 + 1 \otimes t, \qquad
\Delta x = x \otimes 1 + 1 \otimes x - t \otimes v, \qquad
\Delta v = v \otimes 1 + 1 \otimes v; \\[2mm]
 \epsilon(f)=0, \qquad f \in \{t, x, v\}; \\[2mm]
 S(v)= -v,  \quad \quad  S(x)= -x - t v,  \quad \quad  S(t)= -t.
\end{array}
\end{equation}
The nondegenerate pairing between $U_z(\mathfrak{g}(1,1))$ and
$F_z(G(1,1))$ is given by
\begin{equation}
\langle K^m H^n P^p, v^q t^r x^s \rangle =
   m! n! p! \, \, \delta^m_q \delta^n_r \delta^p_s.
\end{equation}

In \cite{olmo00} we constructed the induced representations  of
$U_{z}(\mathfrak{g}(1,1))$, however now we will recover the same results but making use
of its bicrossproduct structure
\begin{equation}
   U_{z}(\mathfrak{g}(1,1)) = U(\mathfrak{v})) \RL F(T_{{z},2}),
\end{equation} 
where $\mathfrak{v}$ is the  Lie algebra of the one-dimensional galilean boosts
group  and $T_{{z}, 2}$ is  a  deformation of the additive group
$\mathbb{R}^2$ defined by
\begin{equation} \label{compodef}
    (b', a')(b,a)= (b' + e^{-2{z} a'}b, a'+a).
\end{equation}
In this definition we have assume that the deformation parameter is
real.
Note that the composition law (\ref{compodef}) is obtained from the  expression
of the coproduct (\ref{estrgalz}).
The elements of $T_{{z}, 2}$  can be factorized as $(b,a)=(b,0)(0,a)$.
The coordinates on  $T_{{z}, 2}$ will be denoted by  $H$ and $P$, so
\begin{equation}
   H(b,a)=b,  \qquad P(b,a)= a .
 \end{equation} 
The   $U(\mathfrak{v}))$--module algebra $F(T_{{z},2})$ is described by the action
 \begin{equation}
       H \ract K= -\frac{1}{4{z}}(1-e^{-4{z} P}), \qquad P \ract K=0. 
\end{equation}
The  vector field associated to this action on
$T_{{z},2}$ is
 \begin{equation} 
\hat{K}= -\frac{1}{4{z}}(1-e^{-4{z} P}) \frac{\partial}{\partial H}. 
\end{equation}

\subsubsect{One-parameter flow}

Let us observe that the vector field $\hat{K}$ has
infinite fixed points  ($(b,0), \ b\in \R$), and $P$ is an invariant.
The  integral curves
 \begin{equation}
  \dot{b}=  -\frac{1}{4{z}}(1-e^{-4{z} a}), \qquad  \dot{a}=  0
\end{equation}
determine the  autonomous system
 \begin{equation}
 b(s)=  -\frac{1}{4{z}}(1-e^{-4{z} c_1}) s + c_2, \qquad
                  a(s)=  c_1. 
\end{equation}
The flow associated to  the vector field $\hat{K}$, deduced from its integral
curves, is
 \begin{equation}
         \Phi^s(b,a)= (b-\frac{1}{4{z}}(1-e^{-4{z} a}) s, a).
 \end{equation}
It is defined for any value of $s$,  giving a global
action of  $\mathfrak{V}$ (the Lie group associated to $\mathfrak{v}$) on $T_{{z},2}$
 \begin{equation}
          e^{sK}\lact(b,a)= (b-\frac{1}{4{z}}(1-e^{-4{z} a}) s, a) .
 \end{equation} 
The  group $T_{{z},2}$ is decomposed in two strata under this action:

\ i) The   set of  points $(b,0)$. Each of them is an orbit with stabilizer the
group $\mathfrak{V}$.
  
ii) The other stratum, constituted by the remaining elements of 
$T_{{z},2}$, is a foliation with sheets 
${\cal O}_a=\{(b, a)|a \in \mathbb{R}^*,\  b \in \mathbb{R}\}$. The   isotopy
group of the point   $(0,a) \in {\cal O}_a$ is $\{ e\}$.
 
\subsubsect{Regular co-spaces} \label{co-spaces2}

Theorem \ref{crbicross} allows to  construct
the  regular co-spaces in a direct and immediately way. Remember that $F_{z}(G(1,1))$ can
be
 described considering elements of  the  form
\begin{equation}
    \phi(b,a),\qquad\qquad \phi \in F(\mathfrak{V}), \quad
       (b,a) \in T_{{z},2},
\end{equation}
instead of the monomial elements
$ v^q t^r x^s$.

For $(F_{z}(G(1,1)), \prec, U_{z}(\mathfrak{g}(1,1)))$ one has 
\begin{equation}\label{rgz1}
 \begin{array}{lll}
(\phi (b,a))\ {\prec}\  e^{s K}&= & \phi (e^{s K} \ \cdot \ ) (b,a), \\[0.2cm] 
(\phi [b,a[) \ {\prec}\ H &= & \phi (b- \frac{1}{4{z}}(1-e^{-4{z} a}) v) (b,a),\\[0.2cm] 
 (\phi (b,a))\ \prec P &= &  \phi a (b,a),
 \end{array}\end{equation} 
and for $(F_{z}(G(1,1)), \succ, U_{z}(\mathfrak{g}(1,1)))$
\begin{equation} \label{rgz2}
 \begin{array}{rll}
  e^{s K} \succ (\phi (b,a))&= & \phi(\, \cdot \, e^{s K})
    (b- \frac{1}{4{z}}(1-e^{-4{z} a}) s ,a), \\[0.2cm] 
  H \succ (\phi (b,a))&= &  b \phi  (b,a), \\[0.2cm] 
  P \succ (\phi (b,a))&=  & a \phi  (b,a).
 \end{array}
\end{equation} 
The elements $(b,a) \in T_{{z},2}$ describe  the  subalgebra of 
$F_{z}(G(1,1))$ generated by $t$ and $x$ which, as in the previous case, is stable under 
the  action  (\ref{rgz2}).

\subsubsect{Induced representations}
A representative of  each   equivalence class of induced representations,
obtained according to the Theorem \ref{repequiv}, is:

\ i) Considering  the  character given by $(b,0)$:
\begin{equation}
   \phi \dashv e^{s K}=  \phi(e^{s K} \, \cdot \,), \qquad
   \phi \dashv H=  \phi b, \qquad
   \phi \dashv P=   0 .
\end{equation}

ii) Taking the  character    associated to $(0,a)$  the induced  representation   
 is
\begin{equation} 
   \phi \dashv e^{s K}=  \phi(e^{s K} \, \cdot \, ), \qquad
   \phi \dashv H=  \phi  \frac{-1}{4{z}}(1-e^{-4{z} a}) v, \qquad
   \phi \dashv P=  \phi a.
 \end{equation}

The local  representations induced  by the  character of 
$U(\mathfrak{so}_0(2))$ given by
\begin{equation} 
K^m \vdash 1= c^m,
\end{equation}
are obtained  applying the   result (\ref{reploc}):
\begin{equation} \begin{split}
     e^{s K} \vdash (b,a)= & \ e^{s c}
                    (b-\frac{1}{4{z}}(1-e^{-4{z} a}) s, a), \\[0.2cm]
     H \vdash (b,a)= &\  b (b,a), \\[0.2cm]
     P \vdash (b,a)= & \ a (b,a).
   \end{split}
\end{equation} 
The actions of  the generators
in the way that they were presented in \cite{olmo00} can be easily deduced from these 
expressions.


\subsect{Quantum kappa--Galilei algebra}
\label{deformacionestandarg}

A contraction of the quantum algebra $U_q(su(2))$ gives the
deformation $U_{\kappa}(\mathfrak{g}(1,1))$  of the enveloping Galilei algebra in $(1+1)$
dimensions \cite{malo98}. 
This quantum algebra is characterized by the following commutation relations and
structure mappings:
\begin{equation}
\begin{array}{c}
[H, K]= - P,  \qquad [P, K]= \frac{P^2}{2\kappa}, \qquad [H, P]=0; \\[2mm]
\Delta H = H \otimes 1 + 1 \otimes H, \qquad
\Delta X = X \otimes 1 + e^{-  H/\kappa} \otimes X,
\quad  X \in \{P, K \}; \\[2mm]
\epsilon(X)= 0, \quad X \in \{H, P, K\}; \\[2mm]
S(H)= -H,  \qquad   S(X)= - e^{{H/\kappa}} X,  \quad  X \in \{P, K\}.
\end{array}
\end{equation}

The bicrossproduct structure of $U_\kappa (\mathfrak{g}(1,1))$ is 
$$ 
U_\kappa (\mathfrak{g}(1,1))={\cal K} \bicross {\cal L} ,
$$
with 
${\cal L}$ the commutative and non-cocommutative Hopf subalgebra  
$U_\kappa (\mathfrak{t}_2)$ spanned by $P$ and $H$, and $\cal K$ the commutative and
cocommutative Hopf subalgebra generated by $K$ (it is not a Hopf  subalgebra of
 $U_\kappa (\mathfrak{g}(1,1))$).
The right action of ${\cal K}$ on ${\cal L}$ is given by
 \begin{equation}
        P \ract K= [P,K]= \frac{P^2}{2\kappa}, \qquad H \ract K= [H,K]= -P ,
 \end{equation}
and the left coaction of ${\cal L}$ over the generator of ${\cal K}$ is
 \begin{equation}
          K \lcact= e^{-  H/\kappa} \otimes K.
 \end{equation}

The dual algebra has also a bicrossproduct structure   
$$
F_\kappa (G(1,1)) ={\cal K}^* \LR {\cal L}^* ,
$$  
where ${\cal K}^*$ is generated by $v$ and ${\cal L}^*$ by  $x$ and $t$.
The left action of  ${\cal L}^*$ on ${\cal K}^*$ is defined by
 \begin{equation}
       x \lact v= \frac{v^2}{2\kappa}, \qquad t \lact v = - v/\kappa,
   \end{equation}
and the right  coaction of ${\cal K}^*$ on ${\cal L}^*$ is:
 \begin{equation}
      \rcact t= t \otimes 1, \qquad \rcact x = x \otimes 1 - t \otimes v.
 \end{equation}
The above action and coaction allow to  recover the Hopf algebra
structure of $F_\kappa (G(1,1))$:
\begin{equation}
\begin{array}{c}
[t, x]= - x/\kappa, \qquad [x,v] = \frac{v^2}{2\kappa}, \qquad [t, v]= - v/\kappa;\\[2mm]
\Delta t = t \otimes 1 + 1 \otimes t, \quad
\Delta x = x \otimes 1 + 1 \otimes x - t \otimes v \quad, 
\Delta v = v \otimes 1 + 1 \otimes v ; \\[2mm]
 \epsilon(f)= 0, \quad f\in \{v, t, x\}; \\[2mm]
 S(v)= -v,  \qquad  S(x)= -x - t v,  \qquad   S(t)= -t.
\end{array}
\end{equation}
The pairing between $U_{\kappa}(\mathfrak{g}(1,1))$ and
$F_{\kappa}(G(1,1))$ is now given by
\begin{equation}
\langle K^m P^n H^p, v^q x^r t^s \rangle =
m! n! p! \,  \delta^m_q \delta^n_r \delta^p_s.
\end{equation}

Let us interpret the algebraic
structure of bicrossproduct of the quantum $\kappa$--Galilei algebra as
\begin{equation}
   U_\kappa(\mathfrak{g}(1,1)) = U(\mathfrak{v}) \RL F(T_{\kappa,2}),
\end{equation} 
where    $T_{\kappa, 2}$ is the   group, deformation of the additive  
$\mathbb{R}^2$ group, defined by
\begin{equation}
    (a', b')(a,b)= (a' + e^{- b'/\kappa}a, b'+b).
\end{equation} 
The elements $(a,b)$ can be factorized in  the  form $(a,b)=(a,0)(0,b)$.
 The functions $P$ and $H$ defined by
\begin{equation}
   P(a,b)=a, \qquad  H(a,b)= b,
\end{equation}
determine a global chart on $T_{\kappa, 2}$.

The  action of the generator $K$ on the  $U(\mathfrak{v})$--module algebra
$F(T_{\kappa,2})$ is given by
 \begin{equation}
H \ract K= - P, \qquad  P \ract K= \frac{P^2}{2\kappa} .
\end{equation} 
Hence, the induced vector field  is
 \begin{equation}
   \hat{K}= \frac{P^2}{2\kappa} \frac{\partial}{\partial P}-
    P   \frac{\partial}{\partial H}.
 \end{equation}


\subsubsect{One-parameter flow} 
 
The invariant points of the vector field $\hat{K}$ are  $(0,b)$. To get an invariant
function under the  action of  $\hat{K}$ is sufficient to determine firstly the  
one--forms, $\eta$,  verifying $\hat{K} \rfloor \eta=0$. The  general solution  is
 \begin{equation}
        \eta_\alpha= \alpha(dP + \frac{1}{2\kappa} P dH),
        \qquad \alpha \in F(T_{\kappa,2}).
 \end{equation}
Choosing $\alpha_0= {1}/{P}$,
the  one--form $\eta_{\alpha_0}$ is exact and invariant under $\hat{K}$.
So,  the  invariant function is
$h= P e^{H/2\kappa}.
$ 
The   autonomous system
 \begin{equation}
  \dot{a}=  \frac{a^2}{2\kappa}, \qquad  \dot{b}=  -a,
 \end{equation}
which determines the integral curves,
is easily integrated. For the curves placed in  the  region $a<0$
we find the following  expressions:
 \begin{equation}
    \begin{split}
    a(s)= & \frac{-1}{c_1 + \frac{s}{2\kappa}}, \\[0.2cm]
    b(s)= & 2\kappa \ln(c_1 + \frac{s}{2\kappa})+ c_2.
    \end{split}
 \end{equation}
The associated second order equations
 \begin{equation} 
    \ddot{a}-\frac{a^3}{2\kappa^2} = 0, \qquad
    \ddot{b}+ \frac{{\dot{b}}^2}{2\kappa}=0,
 \end{equation}
can be interpreted as particles moving in a straight line under  forces depending
on  the  position or  the  velocity, respectively. From  the  expression of  the
 integral curves  the  flow associated to  $\hat{K}$ is obtained:
 \begin{equation}
    \Phi^s(a,b)=
     (\frac{a}{1-\frac{s a}{2\kappa}}, b + 2\kappa \ln(1- \frac{s a}{2\kappa})).
  \end{equation}
Note that  the    action of  $\mathfrak{V}$  on $T_{\kappa,2}$,
 \begin{equation}
     e^{sK}\lact(a,b)=
          (\frac{a}{1-\frac{s a}{2\kappa}}, b + 2\kappa \ln(1- \frac{s a}{2\kappa})),
 \end{equation}
is not global.  The   space $T_{\kappa,2}$ is decomposed in two  strata
under this  action:
 
  \ i) The  set points of  the  form $(0,b)$. Each of them
constitutes a \mbox{$0$--dimensional}  orbit with  stabilizer $\mathfrak{V}$.
  
ii) The   other stratum, constituted by  the   rest of the  space,
        presents a foliation by  one-dimensional sheets.
 
\subsubsect{Regular co-spaces} \label{co-spaces3}
 
 The  action  on the  regular co-space
$(F_\kappa(G(1,1)), \prec, U_\kappa(\mathfrak{g}(1,1)))$
is obtained applying   Theorem  \ref{crbicross}
\begin{equation}
 \begin{split}
   (\phi (a,b))\prec e^{s K}= &\phi (e^{s K} \, \cdot \, ) (a,b), \\[0.2cm]
   (\phi (a,b))\prec P= &  \phi \frac{a}{1- \frac{a v}{2\kappa}} (a,b),\\[0.2cm]
   (\phi (a,b))\prec H= & \phi (b+ 2\kappa \ln(1-\frac{a v}{2\kappa})) (a,b),
 \end{split}
\end{equation}
with $\phi \in F(\mathfrak{V})$ and $(a,b) \in T_{\kappa,2}$. 

The   co-space $(F_\kappa(G(1,1)), \succ, U_\kappa(\mathfrak{g}(1,1)))$
is analogously described by
\begin{equation}
 \begin{split}
  e^{s K} \succ (\phi (a,b))= &\phi(\, \cdot \,  e^{s K})
   (\frac{a}{1-\frac{a s}{2\kappa}}, b + 2\kappa \ln(1- \frac{a s}{2\kappa})), \\[0.2cm]
  P \succ (\phi (a,b))=  & a \phi  (a,b), \\[0.2cm]
  H \succ (\phi (a,b))= &  b \phi  (a,b).
 \end{split}
\end{equation}
Similar comments to those of subsections \ref{co-spaces1} and  \ref{co-spaces2} may be done here.
\subsubsect{Induced representations} 

According to    Theorem \ref{repequiv} each element $(a,b) \in T_{\kappa,2}$ induces a
representation   given by
\begin{equation}
  \begin{split}
     \phi \dashv e^{s K}= & \phi(e^{s K} \; \cdot \;), \\[0.2cm]
     \phi \dashv P= & \phi \frac{a}{1- \frac{a v}{2\kappa}}, \\[0.2cm]
     \phi \dashv H= & \phi (b+ 2\kappa \ln(1-\frac{a v}{2\kappa})),
   \end{split}
\end{equation}
which  effectively  coincides with that was obtained in \cite{olmo00}.

The local  representations induced by  the  character of 
$U(\mathfrak{v})$ given by
\begin{equation} K^m \vdash 1= c^m,\end{equation}
are obtained  applying the   result (\ref{reploc})
\begin{equation} \begin{split}
     e^{s K} \vdash (a,b)= & e^{s c}
 (\frac{a}{1-\frac{a s}{2\kappa}}, b + 2\kappa \ln(1- \frac{a s}{2\kappa})), \\[0.2cm]
     P \vdash (a,b)= & a (b,a), \\[0.2cm]
     H \vdash (a,b)= & b (b,a).
   \end{split}
\end{equation}
From these expressions the  actions of  the generators are easily obtained.

\sect{Concluding remarks}

Remember that in \cite{olmo00} we introduced an algebraic method
for constructing (co)induced representations of Hopf algebras based on the existence of a
triplet composed by two Hopf algebras and a nondegenerate pairing between them such that
there exists a paring of dual bases. However, the difficulty of the computation of the
normal ordering of a product of elements increases with the number of
algebra generators. 
In this work we avoid these troubles when the quantum algebra has a bicrossproduct
structure.

We are able to define structures over a bicrossproduct
Hopf algebra $H= {K}\bicross {L}$ in terms of those of its components ${K}$ and ${L}$.
So, theorem \ref{tbasesduales} gives a procedure to obtain dual bases of the pair ($H,
H^*$) starting from the dual bases of the components. Analogously, theorem
\ref{estructuraestrellabicross}  characterizes a
$*$--structure for the algebra sector of $H$ from the $*$--structures defined on
${K}$ and ${L}$. 

Our induction procedure is not a generalization of the induction method for Lie groups.
We introduce the concept of co-space, which generalizes in an algebraical way the concept of
$G$--space (being $G$ a transformation group),  and we establish the connection between induced
representations and regular co-spaces. There are different procedures for compute regular co-spaces but
the introduction of the endomorphisms associated to the regular actions and the use of adjoint operators
respect to the duality form simplifies extraordinarily the computations \cite{olmo00}.  Note that vector
fields have been used to compute commutators and the advisability of using exponential elements instead
of monomial bases. 
For bicrossproduct Hopf algebras,
like  $H= {K}\bicross {L}$, with $K$ cocommutative and $L$
commutative, theorems \ref{tind0} and \ref{tind1} establish  a connection  between
the representations of $H$  induced by characters of $L$ and certain one-parameter flows.
Although the proof is based on the use of pairs of dual bases the results so obtained are,
essentially, independent of the bases used.   Moreover, we can associate, in some sense, quantum
bicrossproduct groups and dynamical systems via these flows. These relation we will be analyzed more
detailed in a forthcoming paper.

The bicrossproduct Hopf algebras like $H={\cal K}\bicross {\cal L}$, (${\cal K}$ and ${\cal L}$ 
commutative and cocommutative, respectively, infinite dimensional algebras), has been studied
interpreting $\cal K$ as the enveloping algebra $U(\mathfrak k)$ of a Lie group $K$ and $\cal L$ as the
algebra of functions over a Lie group $L$. From this point of view, certain families of Hopf algebras
that are deformations of semidirect products  can be seen as homotopical deformations of the original
actions.

The description of the regular co-spaces associated to $H$  may be done without monomial
bases. Theorem \ref{crbicross0} proves that the action on such co-spaces may be obtained
using the action, deduced from the bicrossproduct structure, of the group $K$ over the
group $L$.  In this way the problems derived from the use of dual bases and the high
dimension of the algebra $H$ are avoided.

The description of the (induced) representations  appears as a corollary of the above
mentioned theorems. The bicrossproduct algebra $H={(\mathfrak k)}\bicross {F(L)}$ gives,
in a natural way, a $*$--structure  for which the representations are, essentially,
unitary. Theorem \ref{repequiv} discusses the equivalence of the induced representations
establishing a correspondence  among classes of induced representations and orbits of $L$
under the action of $K$. This result is in some sense analogous to the Kirillov orbits
method \cite{kirillov}.  The problem of the irreducibility of the representations is still open. Partial
results for particular cases have been obtained; for instance,  see ref. \cite{olmo98} for the
standard quantum $(1+1)$ Galilei algebra and
\cite{bgst98} for the quantum extended $(1+1)$ Galilei algebra. 

\section*{Acknowledgments}
This work has been partially supported by DGES  of the
Ministerio of  Educaci\'{o}n and Cultura of Spain under project PB98-0360, and
the Junta de  Castilla y  Le\'on (Spain).



\begin{thebibliography}{99}

\bibitem{olmo00} O. Arratia and M.A. del Olmo,  {\em J. Math.  Phys.} 
{\bf 41}, (2000) {4817}.

\bibitem{Molnar}   R. Molnar, {\em J. Algebra} {\bf 47}, (1977) 29.

\bibitem{Maj88}  S. Majid,  {\em Class. Quant. Grav.} {\bf 5}, (1988) 1587.

\bibitem{Maj90b}  S. Majid, {\em J. Alg.} {\bf 130}, (1990) 17.

\bibitem{Maj90c} S. Majid, {\em  Pac. J. Math.} {\bf 141}, (1990) 311.

\bibitem{Maj90d}  S. Majid, {\em  Isr. J. Math.} {\bf 72},  (1990) 133.

\bibitem{wigner39}   E. P. Wigner, Ann. Math. {\bf 40}, 149 (1939).

\bibitem{connes} A. Connes, {\em Non-commutative geometry}. Academic
Press, London 1994.

\bibitem{connes01} A. Connes and G. Landi, {\em Commun. Math. Phys.} {\bf 221}, (2001) {141}.

\bibitem{olmo98} O. Arratia and M.A. del
Olmo, {\em Induced representations of quantum groups}.  Anales de F\'{\i}sica,
Monograf\'{\i}as vol. 5 (Ciemat/RSEF, Madrid 1998); math.QA/0110265. 

\bibitem{oscarTesis}  O. Arratia, {\em Induced representations of quantum algebras} Ph.D. Thesis, (in
spanish), (Universidad de Valladolid, Valladolid 1999).

\bibitem{ao00} O. Arratia and M.A. del Olmo, {\em Elements of the theory of induced
representations for quantum groups}, Publicaciones de la RSME vol. 1,
(RSME, Madrid 2000); math.QA/0110266.

\bibitem{dobrev1} V.K.  Dobrev,  {\em  J. Phys. A} {\bf 27}, (1994) 4841 and 6633. 

\bibitem{dobrev2} V.K.  Dobrev,  {\em Induced representations and invariant operators for
quantum groups},  to appear in {\em Turk. J. Phys.}  

\bibitem{ibort} A. Gonz\'alez-Ruiz and L.A. Ibort, {\em  Phys. Lett.}  
{\bf B296}, (1992) 104.

\bibitem{ciccoli98} N. Ciccoli, {\em Induction of quantum group representations},  math/9804138.

\bibitem{maslanka} P. Ma\'slanka,  {\em J. Math. Phys.} {\bf 35}, (1994) 5047.

\bibitem{giller} S. Giller, C. Gonera, P. Kosi\'nski and P.
Ma\'slanka, hep-th 9505007. 

\bibitem{bgst98} F. Bonechi, R. Giachetti, E. Sorace and M. Tarlini,
  {\em Lett. Math. Phys.} {\bf 43}, (1998) 309.

\bibitem{ChP} V. Chari and A. Pressley, {\em A guide to quantum groups}, (Cambridge
Univ. Press, Cambridge 1994).

\bibitem{Maj95a} S. Majid, {\em Foundations of quantum group theory},
(Cambridge Univ. Press, Cambridge 1995).

\bibitem{tmatrix} C. Fronsdal and A. Galindo, {\em Lett. Math. Phys. }  
{\bf 27}, (1993) 39.

\bibitem{Gel43}  I. M. Gelfand and M. A. Naimark,  {\em  Mat. Sbornik} {\bf 12} (1943) 197. 

\bibitem{Bon98}F. Bonechi, N. Ciccoli, R. Giachetti, E. Sorace and M. Tarlini,         
{\em Unitarity of induced representations from coisotropic 
quantum groups},  math/9806062.

\bibitem{olmo84} J.F. Cari\~nena, M.A. del Olmo and M. Santander,  {\em J. Phys. A}  {\bf 17}, (1984)
3091.

\bibitem{Bal95f} A. Ballesteros, E. Celeghini, F. J. Herranz, M. A. del
Olmo  and M. Santander,   {\em J. Phys. A} {\bf 28}, (1995) 3129.

\bibitem{Bal95g}  A. Ballesteros,
F. J. Herranz, C. M. Pere\~na, M. A. del Olmo and  M. Santander,  {\em J. Phys. A} {\bf 28}, 
 (1995) 7113.

\bibitem{PerJ96} J. C. P\'erez Bueno, {\em Contractions of deformed Hopf algebras and
their structure}, M.S. Thesis, (in spanish), (Universidad de 
Valencia 1998).

\bibitem{cgst92} E. Celeghini, R. Giachetti, E. Sorace and M. Tarlini,
{\em Contractions of quantum groups}. Lecture Notes in Mathematics vol. 1510,
 pg. 221 (Springer, Berlin 1992).

\bibitem{cgst91} E. Celeghini, R. Giachetti, E. Sorace and M. Tarlini,
 {\em J. Math. Phys.} {\bf 32} (1991), 1155.

\bibitem{apb96} J. A. of  Azc\'arraga and J. C. P\'erez Bueno,  {\em J. Phys. A} {\bf 29}
(1996), 6353.

\bibitem{malo98} J. C. P\'erez Bueno, {\em Cohomological aspects of some physical
problems}. Ph.D. Thesis (Universidad de  Valencia 1998).

\bibitem{kirillov} A.A. Kirillov, {\sl Elements of the Theory of Representations}
(Springer, Berlin 1976).

\end{thebibliography}
\end{document}